\documentclass[a4paper,12pt]{amsart}
\usepackage[pagebackref=true]{hyperref}
\usepackage{amssymb,amsmath,graphicx,stmaryrd,fullpage,setspace,microtype}
\usepackage[pagebackref=true]{hyperref}
\hypersetup{backref}
\usepackage[all]{xy}
\title{The choice of cofibrations of higher dimensional
  transition systems}

\author[P. Gaucher]{Philippe Gaucher}

\address{Laboratoire PPS  (CNRS UMR 7126)\\
  Univ Paris Diderot\\
Sorbonne Paris Cit\'e\\
 Case 7014\\ 75205 PARIS Cedex 13 \\ France}

\urladdr{http://www.pps.univ-paris-diderot.fr/{\~{}}gaucher/} 

\subjclass{18C35,18G55,55U35,68Q85}

\keywords{left determined model category, combinatorial model category, higher dimensional transition system, causal structure}

\swapnumbers

\newcommand{\C}{\mathcal{C}}

\newcommand{\K}{\mathcal{K}}
\newcommand{\LL}{\mathcal{A}}
\newcommand{\I}{\mathcal{I}}

\newcommand{\de}{\partial}
\newcommand{\p}\times

\newtheorem{thm}{Theorem}[section]
\newtheorem{prop}[thm]{Proposition}
\newtheorem{lem}[thm]{Lemma}
\newtheorem{cor}[thm]{Corollary}

\newtheorem{defn}[thm]{Definition}
\newtheorem{nota}[thm]{Notation}

\newcommand{\bd}{\begin{defn}}
\newcommand{\ed}{\end{defn}}
\newcommand{\bp}{\begin{prop}}
\newcommand{\ep}{\end{prop}}
\newcommand{\bth}{\begin{thm}}
\renewcommand{\eth}{\end{thm}}
\newcommand{\bpf}{\begin{proof}}
\newcommand{\epf}{\end{proof}}
\newcommand{\bc}{\begin{cor}}
\newcommand{\ec}{\end{cor}}

\newcommand{\fL}[1]{\ar@{->}[ll]_-{#1}}
\newcommand{\fR}[1]{\ar@{->}[rr]^-{#1}}
\newcommand{\fRr}[1]{\ar@{->}[rrr]^-{#1}}
\newcommand{\fD}[1]{\ar@{->}[dd]_-{#1}}
\newcommand{\fU}[1]{\ar@{->}[uu]^-{#1}}
\newcommand{\f}[2]{\ar@{->}[#1]|{#2}}
\newcommand{\ff}[2]{\ar@2{->}[#1]|{#2}}
\newcommand{\frr}[1]{\ar@{->}[rrrr]^-{#1}}

\newcommand{\fl}[1]{\ar@{->}[l]_-{#1}}
\newcommand{\fr}[1]{\ar@{->}[r]^-{#1}}
\newcommand{\fd}[1]{\ar@{->}[d]_-{#1}}
\newcommand{\fu}[1]{\ar@{->}[u]^-{#1}}

\newcommand{\iso}{\cong}

\renewcommand{\leq}{\leqslant}
\renewcommand{\geq}{\geqslant}
\newcommand{\dd}[1]{\uparrow\!\!{#1}\!\!\uparrow}

\def\cartesien{%
  \ar@{-}[]+R+<6pt,-2pt>;[]+RD+<6pt,-6pt>%
  \ar@{-}[]+D+<2pt,-6pt>;[]+RD+<6pt,-6pt>%
}
\def\cocartesien{%
  \ar@{-}[]+L+<-6pt,+2pt>;[]+LU+<-6pt,+6pt>%
  \ar@{-}[]+U+<-2pt,+6pt>;[]+LU+<-6pt,+6pt>%
}
\def\hocartesien{%
  \ar@{-}[]+R+<6pt,-2pt>;[]+RD+<6pt,-6pt>_{h}%
  \ar@{-}[]+D+<2pt,-6pt>;[]+RD+<6pt,-6pt>%
}
\def\hococartesien{%
  \ar@{-}[]+L+<-6pt,+2pt>;[]+LU+<-6pt,+6pt>_{h}%
  \ar@{-}[]+U+<-2pt,+6pt>;[]+LU+<-6pt,+6pt>%
}

\newcommand{\brm}[1]{\rm{\mathbf{#1}}}

\newcommand{\set}{{\brm{Set}}}

\DeclareMathOperator{\id}{Id}

\DeclareMathOperator{\Mor}{Mor}

\newcommand{\liminj}{\varinjlim}

\newcommand{\wts}{\mathcal{W\!T\!S}}
\newcommand{\cts}{\mathcal{C\!T\!S}}

\newcommand{\rts}{\mathcal{R\!T\!S}}
\newcommand{\restr}{\!\upharpoonright\!}

\DeclareMathOperator{\dom}{dom}

\makeatletter
\def\varholim@#1#2{%
  \vtop{\m@th\ialign{##\cr
    \hfil$#1\operator@font holim$\hfil\cr
    \noalign{\nointerlineskip\kern1.5\ex@}#2\cr
    \noalign{\nointerlineskip\kern-\ex@}\cr}}%
}
\def\holimproj{%
  \mathop{\mathpalette\varholim@{\leftarrowfill@\textstyle}}\nmlimits@
}
\def\holiminj{%
  \mathop{\mathpalette\varholim@{\rightarrowfill@\textstyle}}\nmlimits@
}
\makeatother

\DeclareMathOperator{\cell}{{\brm{cell}}}
\DeclareMathOperator{\cof}{{\brm{cof}}}

\DeclareMathOperator{\inj}{{\brm{inj}}}
\newcommand{\ddownarrow}{{\downarrow}}

\DeclareMathOperator{\cyl}{{Cyl}}
\DeclareMathOperator{\cocyl}{{Path}}
\DeclareMathOperator{\CSA}{CSA}

\begin{document}

\begin{abstract} 
  It is proved that there exists a left determined model structure of
  weak transition systems with respect to the class of monomorphisms
  and that it restricts to left determined model structures on cubical
  and regular transition systems. Then it is proved that, in these
  three model structures, for any higher dimensional transition system
  containing at least one transition, the fibrant replacement contains
  a transition between each pair of states. This means that the
  fibrant replacement functor does not preserve the causal
  structure. As a conclusion, we explain why working with star-shaped
  transition systems is a solution to this problem.
\end{abstract}

\maketitle
\tableofcontents

\section{Introduction}

\subsection{Summary of the paper}
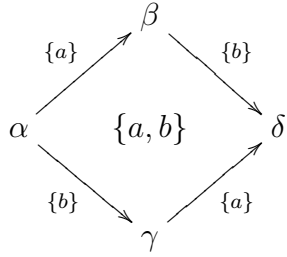
\begin{figure}
\[
\xymatrix{
& \beta \ar@{->}[rd]^{\{b\}}&\\
\alpha\ar@{->}[ru]^{\{a\}}\ar@{->}[rd]_{\{b\}} & \{a,b\} & \delta\\
&\gamma\ar@{->}[ru]_{\{a\}}&}
\]
\caption{$a|| b$ : Concurrent execution of $a$ and $b$}
\label{concab}
\end{figure}

This work belongs to our series of papers devoted to \emph{higher
  dimensional transition systems}. It is a (long) work in progress.
The notion of higher dimensional transition system dates back to
Cattani-Sassone's paper \cite{MR1461821}. These objects are a higher
dimensional analogue of the computer-scientific notion of labelled
transition system. Their purpose is to model the concurrent execution
of $n$ actions by a multiset of actions, i.e. a set with a possible
repetition of some elements (e.g. $\{0,0,2,3,3,3\}$).  The higher
dimensional transition $a||b$ modeling the concurrent execution of the
two actions $a$ and $b$, depicted by Figure~\ref{concab}, consists of
the transitions $(\alpha,\{a\},\beta)$, $(\beta,\{b\},\delta)$,
$(\alpha,\{b\},\gamma)$, $(\gamma,\{a\},\delta)$ and
$(\alpha,\{a,b\},\delta)$. The labelling map is the identity map.
Note that with $a=b$, we would get the $2$-dimensional transition
$(\alpha,\{a,a\},\delta)$ which is not equal to the $1$-dimensional
transition $(\alpha,\{a\},\delta)$. The latter actually does not exist
in Figure~\ref{concab}. Indeed, the only $1$-dimensional transitions
labelled by the multiset $\{a\}$ are $(\alpha,\{a\},\beta)$ and
$(\gamma,\{a\},\delta)$.  The new formulation introduced in
\cite{hdts} enabled us to interpret them as a small-orthogonality
class of a locally finitely presentable category $\wts$ of \emph{weak
  transition systems} equipped with a topological functor towards a
power of the category of sets. In this new setting, the
$2$-dimensional transition of Figure~\ref{concab} becomes the tuple
$(\alpha,a,b,\delta)$. The set of transitions has therefore to satisfy
the multiset axiom (here: if the tuple $(\alpha,a,b,\delta)$ is a
transition, then the tuple $(\alpha,b,a,\delta)$ has to be a
transition as well) and the patching axiom which is a topological
version (in the sense of topological functors) of Cattani-Sassone's
interleaving axiom. We were then able to state a categorical
comparison theorem between them and (labelled) symmetric precubical
sets in \cite{hdts}. We studied in \cite{cubicalhdts} a homotopy
theory of \emph{cubical transition systems} $\cts$ and in \cite{csts},
exhaustively, a homotopy theory of \emph{regular transition systems}
$\rts$. The adjective cubical means that the weak transition system is
the union of its subcubes. In particular this means that every higher
dimensional transition has lower dimensional faces. However, a square
for example may still have more than four $1$-dimensional faces in the
category of cubical transition systems. A cubical transition system is
by definition regular if every higher dimensional transition has the
expected number of faces. All known examples coming from process
algebra are cubical because they are colimits of cubes, and therefore
are equal to the union of their subcubes. Indeed, the associated
higher dimensional transition systems are realizations in the sense of
\cite[Theorem~9.2]{hdts} (see also
\cite[Theorem~7.4]{homotopyprecubical}) of a labelled precubical set
obtained by following the semantics expounded in \cite{ccsprecub}.  It
turns out that there exist colimits of cubes which are not regular
(see the end of \cite[Section~2]{csts}). However, it can also be
proved that all process algebras for any synchronization algebra give
rise to regular transition systems. The regular transition systems
seem to be the only interesting ones. However, their mathematical
study requires to use the whole chain of inclusion functors $\rts
\subset \cts \subset \wts$.

The homotopy theories studied in \cite{cubicalhdts} and in \cite{csts}
are obtained by starting from a left determined model structure on
weak transition systems with respect to the class of maps of weak
transition systems which are one-to-one on the set of actions (but not
necessarily one-to-one on the set of states) and then by restricting
it to full subcategories (the coreflective subcategory of cubical
transition systems, and then the reflective subcategory of regular
ones).

In this paper, we will start from the left determined model category
of weak transition systems with respect to the class of monomorphisms
of weak transition systems, i.e. the cofibrations are one-to-one not
only on the set of actions, but also on the set of states. Indeed, it
turns out that such a model structure exists: it is the first result
of this paper (Theorem~\ref{standard-wts}). And it turns out that it
restricts to the full subcategories of cubical and regular transition
systems as well and that it gives rise to two new left determined
model structures: it is the second result of this paper
(Theorem~\ref{standard-cts} for cubical transition systems and
Theorem~\ref{standard-rts} for regular transition systems).

Unlike the homotopy structures studied in \cite{cubicalhdts} and in
\cite{csts}, the model structures of this paper do not have the map
$R:\{0,1\} \to \{0\}$ identifying two states as a cofibration
anymore. However, there are still cofibrations of regular transition
systems which identify two different states. This is due to the fact
that the set of states of a colimit of regular transition systems is
in general not the colimit of the sets of states. There are
identifications inside the set of states which are forced by the
axioms satisfied by regular transition systems, actually CSA2. This
implies that the class of cofibrations of this new left determined
model structure on regular transition systems, like the one described
and studied in \cite{csts}, still contains cofibrations which are not
monic: see an example at the very end of
Section~\ref{restriction-modelstructure}.

Without additional constructions, these new model structures are
irrelevant for concurrency theory. Indeed, the fibrant replacement
functor, in any of these model categories (the weak transition systems
and also the cubical and the regular ones), destroys the causal
structure of the higher dimensional transition system: this is the
third result of this paper (Theorem~\ref{path-wts} and
Theorem~\ref{path-cts}).

We open this new line of research anyway because of the following
discovery: by working with star-shaped transition systems, the bad
behavior of the fibrant replacement just disappears. This point is
discussed in the very last section of the paper. The fourth result of
this paper is that left determined model structures can be constructed
on star-shaped (weak or cubical or regular) transition systems
(Theorem~\ref{star-shaped}).  This paper is the starting point of the
study of these new homotopy theories.

Appendix~\ref{relocation} is a technical tool to relocate the map
$R:\{0,1\} \to \{0\}$ in a transfinite composition. Even if this map
is not a cofibration in this paper, it still plays an important role
in the proofs. This map seems to play an ubiquitous role in our
homotopy theories.

\subsection{Prerequisites and notations}

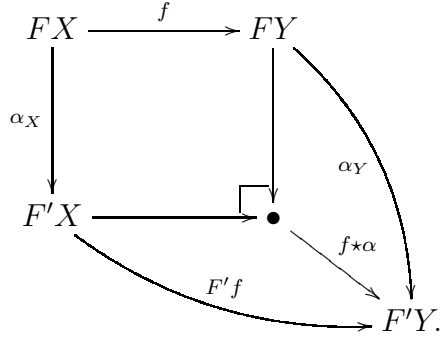
\begin{figure}
\[
\xymatrix{
FX \fD{\alpha_X}\fR{f} && FY \fD{}\ar@/^15pt/@{->}[dddr]_-{\alpha_Y} &\\
&& &&\\
F'X \ar@/_15pt/@{->}[rrrd]^-{F'f}\fR{} &&  \bullet \cocartesien\ar@{->}[rd]^-{f\star \alpha} & \\
&& & F'Y.
}
\]
\caption{Definition of $f\star \alpha$.}
\label{star}
\end{figure}

All categories are locally small. The set of maps in a category $\K$
from $X$ to $Y$ is denoted by $\K(X,Y)$. The class of maps of a
category $\K$ is denoted by $\Mor(\K)$. The composite of two maps is
denoted by $fg$ instead of $f \circ g$. The initial (final resp.)
object, if it exists, is always denoted by $\varnothing$ ($\mathbf{1}$
resp.). The identity of an object $X$ is denoted by $\id_X$.  A
subcategory is always isomorphism-closed.  Let $f$ and $g$ be two maps
of a locally presentable category $\K$. Write $f\square g$ when $f$
satisfies the \emph{left lifting property} (LLP) with respect to $g$,
or equivalently $g$ satisfies the \emph{right lifting property} (RLP)
with respect to $f$. Let us introduce the notations $\inj_\K(\C) = \{g
\in \K, \forall f \in \C, f\square g\}$ and $\cof_\K(\C) = \{f \in \K,
\forall g\in \inj_\K(\C), f\square g\}$ where $\C$ is a class of maps
of $\K$. The class of morphisms of $\K$ that are transfinite
compositions of pushouts of elements of $\C$ is denoted by
$\cell_\K(\C)$. There is the inclusion $\cell_\K(K)\subset
\cof_\K(K)$. Moreover, every morphism of $\cof_\K(K)$ is a retract of
a morphism of $\cell_\K(K)$ as soon as the domains of $K$ are small
relative to $\cell_\K(K)$ \cite[Corollary~2.1.15]{MR99h:55031},
e.g. when $\K$ is locally presentable.  A class of maps of $\K$ is
cofibrantly generated if it is of the form $\cof_\K(S)$ for some set
$S$ of maps of $\K$. For every map $f:X \to Y$ and every natural
transformation $\alpha : F \to F'$ between two endofunctors of $\K$,
the map $f\star \alpha$ is defined by the diagram of
Figure~\ref{star}.  For a set of morphisms $\mathcal{A}$, let
$\mathcal{A} \star \alpha = \{f\star \alpha, f\in \mathcal{A}\}$.

We refer to \cite{MR95j:18001} for locally presentable categories, to
\cite{MR2506258} for combinatorial model categories, and to
\cite{topologicalcat} for topological categories, i.e. categories
equipped with a topological functor towards a power of the category of
sets.  We refer to \cite{MR99h:55031} and to \cite{ref_model2} for
model categories. For general facts about weak factorization systems,
see also \cite{ideeloc}. The reading of the first part of
\cite{MOPHD}, published in \cite{MO}, is recommended for any reference
about good, cartesian, and very good cylinders.

We use the paper \cite{leftdet} as a toolbox for constructing the
model structures. To keep this paper short, we refer to \cite{leftdet}
for all notions related to Olschok model categories.

\section{The model structure of weak transition systems}

We are going first to recall a few facts about weak transition
systems.

\begin{nota} Let $\Sigma$ be a fixed nonempty set of {\rm
    labels}.  \end{nota}

\bd\label{def_HDTS} A {\rm weak transition system} consists of a
triple \[X=(S,\mu:L\rightarrow \Sigma,T=\bigcup_{n\geq 1}T_n)\] where
$S$ is a set of {\rm states}, where $L$ is a set of {\rm actions},
where $\mu:L\rightarrow \Sigma$ is a set map called the {\rm labelling
  map}, and finally where $T_n\subset S\p L^n\p S$ for $n \geq 1$ is a
set of {\rm $n$-transitions} or {\rm $n$-dimensional transitions} such
that the two following axioms hold:
\begin{itemize}
\item (Multiset axiom) For every permutation $\sigma$ of
  $\{1,\dots,n\}$ with $n\geq 2$, if the tuple
  $(\alpha,u_1,\dots,u_n,\beta)$ is a transition, then the tuple
  $(\alpha,u_{\sigma(1)}, \dots, u_{\sigma(n)}, \beta)$ is a
  transition as well.
\item (Patching axiom~\footnote{This axiom is called the Coherence
    axiom in \cite{hdts} and \cite{cubicalhdts}, and the composition
    axiom in \cite{csts}. I definitively adopted the terminology
    ``patching axiom'' after reading the Web page in nLab devoted to
    higher dimensional transition systems and written by Tim Porter.})
  For every $(n+2)$-tuple $(\alpha,u_1,\dots,u_n,\beta)$ with $n\geq
  3$, for every $p,q\geq 1$ with $p+q<n$, if the five tuples
\begin{align*}
 &(\alpha,u_1, \dots, u_n, \beta), \\ &(\alpha,u_1, \dots, u_p,
  \nu_1), (\nu_1, u_{p+1}, \dots, u_n, \beta),\\ &(\alpha, u_1,
  \dots, u_{p+q}, \nu_2), (\nu_2, u_{p+q+1}, \dots, u_n, \beta)
\end{align*}
  are transitions, then the $(q+2)$-tuple $(\nu_1, u_{p+1}, \dots,
  u_{p+q}, \nu_2)$ is a transition as well.
\end{itemize}
A map of weak transition systems
\[f:(S,\mu : L \rightarrow \Sigma,(T_n)_{n\geq 1}) \rightarrow
(S',\mu' : L' \rightarrow \Sigma ,(T'_n)_{n\geq 1})\] consists of a
set map $f_0: S \rightarrow S'$, a commutative square
\[
\xymatrix{
  L \ar@{->}[r]^-{\mu} \ar@{->}[d]_-{\widetilde{f}}& \Sigma \ar@{=}[d]\\
  L' \ar@{->}[r]_-{\mu'} & \Sigma}
\] 
such that if $(\alpha,u_1,\dots,u_n,\beta)$ is a transition, then
$(f_0(\alpha),\widetilde{f}(u_1),\dots,\widetilde{f}(u_n),f_0(\beta))$
is a transition. The corresponding category is denoted by $\wts$.  The
$n$-transition $(\alpha,u_1,\dots,u_n,\beta)$ is also called a {\rm
  transition from $\alpha$ to $\beta$}: $\alpha$ is the initial state
and $\beta$ the final state of the transition.  The maps $f_0$ and
$\widetilde{f}$ are sometimes denoted simply as $f$. \ed

The category $\wts$ is locally finitely presentable and the
functor \[\omega : \wts \longrightarrow \set^{\{s\}\cup \Sigma},\]
where $s$ is the sort of states, taking the weak higher dimensional
transition system $(S,\mu : L \rightarrow \Sigma,(T_n)_{n\geq 1})$ to
the $(\{s\}\cup \Sigma)$-tuple of sets $(S,(\mu^{-1}(x))_{x\in
  \Sigma}) \in \set^{\{s\}\cup \Sigma}$ is topological by
\cite[Theorem~3.4]{hdts}. The terminal object of $\wts$ is the weak
transition system
\[\mathbf{1}=(\{0\},\id_\Sigma:\Sigma\to\Sigma,\bigcup_{n\geq 1}\{0\}
\p \Sigma^n \p \{0\}).\]

\begin{nota} For $n\geq 1$, let $0_n = (0,\dots,0)$ ($n$ times) and
  $1_n = (1,\dots,1)$ ($n$ times). By convention, let
  $0_0=1_0=()$. \end{nota}

Here are some important examples of weak transition systems: 
\begin{enumerate}
\item Every set $S$ can be identified with the weak transition system
  having the set of states $S$, with no actions and no
  transitions. For all weak transition system $X$, the set
  $\wts(\{0\},X)$ is the set of states of $X$. The empty set is the
  initial object of $\wts$.
\item The weak transition system $\underline{x} = (\varnothing, \{x\}
  \subset \Sigma, \varnothing)$ for $x\in \Sigma$. For all weak
  transition system $X$, the set $\wts(\underline{x},X)$ is the set of
  actions of $X$ labelled by $x$ and $\bigsqcup_{x\in
    \Sigma}\wts(\underline{x},X)$ is the set of actions of $X$.
\item Let $n \geq 0$. Let $x_1,\dots,x_n \in \Sigma$. The \emph{pure
    $n$-transition} $C_n[x_1,\dots,x_n]^{ext}$ is the weak transition
  system with the set of states $\{0_n,1_n\}$, with the set of actions
  \[\{(x_1,1), \dots, (x_n,n)\}\] and with the transitions all
  $(n+2)$-tuples \[(0_n,(x_{\sigma(1)},\sigma(1)), \dots,
  (x_{\sigma(n)},\sigma(n)),1_n)\] for $\sigma$ running over the set of
  permutations of the set $\{1,\dots ,n\}$. Intuitively, the pure
  transition is a cube without faces of lower dimension.  For all weak
  transition system $X$, the set $\wts(C_n[x_1,\dots,x_n]^{ext},X)$ is
  the set of transitions $(\alpha,u_1,\dots,u_n,\beta)$ of $X$ such
  that for all $1\leq i\leq n$, $\mu(u_i)=x_i$ and
  \[\bigsqcup_{x_1,\dots,x_n\in \Sigma}\wts(C_n[x_1,\dots,x_n]^{ext},X)\]
  is the set of transitions of $X$.
\end{enumerate}

The purpose of this section is to prove the existence of a left
determined combinatorial model structure on the category of weak
transition systems with respect to the class of monomorphisms. 

We first have to check that the class of monomorphisms of weak
transition systems is generated by a set. The set of generating
cofibrations is obtained by removing the map $R:\{0,1\}\to \{0\}$ from
the set of generating cofibrations of the model structure studied in
\cite{cubicalhdts} and in \cite{csts}.

\begin{nota} (Compare with
  \cite[Notation~5.3]{cubicalhdts}) \label{defI} Let $\I$ be the set
  of maps $C:\varnothing \rightarrow \{0\}$, $\varnothing \subset
  \underline{x}$ for $x\in \Sigma$ and $\{0_n,1_n\} \sqcup
  \underline{x_1} \sqcup \dots \sqcup \underline{x_n} \subset C_n[x_1,
  \dots, x_n]^{ext}$ for $n\geq 1$ and $x_1, \dots, x_n \in
  \Sigma$. \end{nota}

\begin{lem} \label{colim-preserv0} The forgetful functor mapping a
  weak transition system to its set of states is colimit-preserving.
  The forgetful functor mapping a weak transition system to its set of
  actions is colimit-preserving.
\end{lem}

\bpf The lemma is a consequence of the fact that the forgetful functor
$\omega : \wts \longrightarrow \set^{\{s\}\cup \Sigma}$ taking the
weak higher dimensional transition system $(S,\mu : L \rightarrow
\Sigma,(T_n)_{n\geq 1})$ to the $(\{s\}\cup \Sigma)$-tuple of sets
$(S,(\mu^{-1}(x))_{x\in \Sigma}) \in \set^{\{s\}\cup \Sigma}$ is
topological.  \epf

\begin{lem} \label{R-wts-epic}
All maps of $\cell_\wts(\{R\})$ are epic.
\end{lem}

\bpf Let $f,g,h$ be three maps of $\wts$ with $f\in \cell_\wts(\{R\})$
such that $gf=hf$. By functoriality, we obtain the equality
$\omega(g)\omega(f)=\omega(h)\omega(f)$.  All maps of
$\cell_\wts(\{R\})$ are onto on states and the identity on actions by
Lemma~\ref{colim-preserv0}. Therefore $\omega(f)$ is epic and we
obtain $\omega(g)=\omega(h)$. Since the forgetful functor $\omega :
\wts \longrightarrow \set^{\{s\}\cup \Sigma}$ is topological, it is
faithful by \cite[Theorem~21.3]{topologicalcat}. Thus, we obtain
$g=h$.  \epf

\bp\label{small-mono} There is the equality $\cell_{\wts}(\I) =
\cof_{\wts}(\I)$ and this class of maps is the class of monomorphisms
of weak transition systems.  \ep

\bpf By \cite[Proposition~3.1]{cubicalhdts}, a map of weak transition
systems is a mono\-morphism if and only if it induces a one-to-one set
map on states and on actions. Consequently, by
\cite[Proposition~5.4]{cubicalhdts}, a cofibration of weak transition
systems $f$ belongs to $\cell_{\wts}(\I \cup \{R\})$. All maps of $\I$
belong to $\inj_\wts(\{R\})$ because they are one-to-one on states.
Using Lemma~\ref{R-wts-epic}, we apply Theorem~\ref{move_bad}: $f$
factors uniquely, up to isomorphism, as a composite $f=f^+f^-$ with
$f^+\in \cell_{\wts}(\I)$ and $f^-\in \cell_{\wts}(\{R\})$.  The map
$f^-$ is one-to-one on states because $f$ is one-to-one on states. We
obtain the equalities $f^-=\id$ and $f=f^+$. Therefore $f$ belongs to
$\cell_{\wts}(\I)$.  Conversely, every map of $\cell_{\wts}(\I)$ is
one-to-one on states and on actions by
Lemma~\ref{colim-preserv0}. Thus, the class of cofibrations is
$\cell_{\wts}(\I)$. Since the underlying category $\wts$ is locally
presentable, every map of $\cof_{\wts}(\I)$ is a retract of a map of
$\cell_{\wts}(\I)$. This implies that every map of $\cof_{\wts}(\I)$
is one-to-one on states and actions. Thus, we obtain $\cof_{\wts}(\I)
\subset \cell_{\cts}(\I)$. Hence we have obtained $\cof_{\wts}(\I) =
\cell_{\wts}(\I)$ and the proof is complete.  \epf

Let us now introduce the interval object of this model structure.

\bd Let $V$ be the weak transition system defined as
follows:
\begin{itemize}
\item The set of states is $\{0,1\}$.
\item The set of actions is $\Sigma\p \{0,1\}$.
\item The labelling map is the projection $\Sigma\p \{0,1\} \to \Sigma$.
\item The transitions are the tuples
  \[(\epsilon_0,(x_1,\epsilon_1),\dots,(x_n,\epsilon_n),\epsilon_{n+1})\]
  for all $\epsilon_0,\dots,\epsilon_{n+1}\in \{0,1\}$ and all $x_1,\dots,x_n\in \Sigma$.
\end{itemize}
\ed

\begin{nota} Denote by $\cyl :\wts \to \wts$ the {\rm functor} $-\p
  V$.
\end{nota}

\bp \label{cyl-construction} Let $X=(S,\mu:L\to \Sigma,T)$ be a weak
transition system. The weak transition system $\cyl(X)$ has the set of
states $S\p \{0,1\}$, the set of actions $L\p \{0,1\}$, the labelling
map the composite map $\mu:L\p \{0,1\} \to L \to \Sigma$, and a tuple
\[((\alpha,\epsilon_0),(u_1,\epsilon_1),\dots,(u_n,\epsilon_n),(\beta,\epsilon_{n+1}))\]
is a transition of $\cyl(X)$ if and only if the tuple
$(\alpha,u_1,\dots,u_n,\beta)$ is a transition of $X$.  There exists a
unique map of weak transition systems $\gamma_X^\epsilon:X\to\cyl(X)$
for $\epsilon=0,1$ defined on states by $s\mapsto (s,\epsilon)$ and on
actions by $u\mapsto (u,\epsilon)$.  There exists a unique map of weak
transition systems $\sigma_X:\cyl(X)\to X$ defined on states by
$(s,\epsilon)\mapsto s$ and on actions by $(u,\epsilon)\mapsto u$.
There is the equality $\sigma_X \gamma_X^\epsilon = \id_X$. The
composite map $\sigma_X \gamma_X$ with $\gamma_X=\gamma_X^0 \sqcup
\gamma_X^1$ is the codiagonal of $X$.  \ep

Note that if $T_n$ denotes the set of $n$-transitions of $X$, then the
set of $n$-transitions of $\cyl(X)$ is $T_n \p \{0,1\}^{n+2}$.

\bpf The binary product in $\wts$ is described in
\cite[Proposition~5.5]{cubicalhdts}. The set of states of $\cyl(X)$ is
$S\p \{0,1\}$. The set of actions of $\cyl(X)$ is the product $L
\p_\Sigma (\Sigma \p \{0,1\}) \iso L \p \{0,1\}$ and the transitions
of $\cyl(X)$ are the tuples of the form
$((\alpha,\epsilon_0),(u_1,\epsilon_1),\dots,(u_n,\epsilon_n),(\beta,\epsilon_{n+1}))$
such that $(\alpha,u_1,\dots,u_n,\beta)$ is a transition of $X$ and
such that the tuple
$(\epsilon_0,(\mu(u_1),\epsilon_1),\dots,(\mu(u_n),\epsilon_n),\epsilon_{n+1})$
is a transition of $V$. The latter holds for any choice of
$\epsilon_0,\dots,\epsilon_{n+1}\in \{0,1\}$ by definition of $V$.
\epf

\bp Let $X$ be a weak transition system. Then the map
$\gamma_X:X\sqcup X\to \cyl(X)$ is a monomorphism of weak
transition systems and the map $\sigma_X:\cyl(X) \to X$
satisfies the right lifting property (RLP) with respect to the
monomorphisms of weak transition systems. \ep

\bpf By \cite[Proposition~3.1]{cubicalhdts}, the map $\gamma_X:X\sqcup
X\to \cyl(X)$ is a monomorphism of $\wts$ since it is
bijective on states and on actions.  The lift $\ell$ exists in the
following diagram:
\[
\xymatrix
{
\varnothing \fR{}\fD{C} && {V} \fD{} \\
&& \\
\{0\} \fR{} \ar@{-->}[rruu]^-{\ell}&& \mathbf{1}
}
\]
where $\mathbf{1}=(\{0\},\id_\Sigma:\Sigma\to\Sigma,\bigcup_{n\geq
  1}\{0\} \p \Sigma^n \p \{0\})$ is the terminal object of $\wts$:
take $\ell(0)=0$.  The lift $\ell$ exists in the following diagram:
\[
\xymatrix
{
\varnothing \fR{}\fD{} && {V} \fD{} \\
&& \\
\underline{x} \fR{} \ar@{-->}[rruu]^-{\ell}&& \mathbf{1}.
}
\]
Indeed, $\ell(x)=x$ is a solution.
Finally, consider a commutative diagram of the form: 
\[
\xymatrix
{
\{0_n,1_n\} \sqcup \underline{x_1} \sqcup \dots \sqcup
  \underline{x_n}  \fR{\phi}\fD{\subset} && {V} \fD{} \\
&& \\
C_n[x_1, \dots, x_n]^{ext} \fR{} \ar@{-->}[rruu]^-{\ell}&& \mathbf{1}.
}
\]
Then let
\[\ell(0_n,(x_{\sigma(1)},1),\dots,(x_{\sigma(n)},n),1_n)=(\phi(0_n),(x_{\sigma(1)},0),\dots,(x_{\sigma(n)},0),\phi(1_n))\]
for any permutation $\sigma$: it is a solution. Therefore by
Proposition~\ref{small-mono}, the map ${V}\to\mathbf{1}$ satisfies the
RLP with respect to all monomorphisms. Finally, consider the
commutative diagram of solid arrows:
\[
\xymatrix
{
A \fD{f} \fR{} && \cyl(X) \fD{\sigma_X} \\
&& \\
B \ar@{-->}[rruu]^-{\ell} \fR{} && X
}
\]
where $f$ is a monomorphism. Then the lift $\ell$ exists because there
are the isomorphisms $\cyl(X)\iso X\p V$ and
$X\iso X\p \mathbf{1}$ and because the map $\sigma_X$ is equal to the
product $\id_X\p (V\to \mathbf{1})$.  \epf

\begin{cor}
  The functor $\cyl:\wts \to \wts$ together with the natural
  transformations $\gamma:\id \Rightarrow \cyl$ and $\cyl\Rightarrow
  \id$ gives rise to a very good cylinder with respect to $\I$.
\end{cor}

\bp \label{cyl-colim} The functor $\cyl:\wts \to \wts$ is
colimit-preserving. \ep

We will use the following notation: let $I$ be a small category. For
any diagram $D$ of weak transition systems over $I$, the canonical map
$D_i \to \liminj D_i$ is denoted by $\phi_{D,i}$.

\bpf Let $I$ be a small category. Let $X:i\mapsto X_i$ be a small
diagram of weak transition systems over $I$. By
Lemma~\ref{colim-preserv0}, for all objects $i$ of $I$, the map
$\phi_{X,i}:X_i \to \liminj_i X_i$ is the inclusion $S_i \subset
\liminj_i S_i$ on states and the inclusion $L_i \subset \liminj_i L_i$
on actions if $S_i$ ($L_i$ resp.)  is the set of states (of actions
resp.) of $X_i$. By definition of the functor $\cyl$, for all objects
$i$ of $I$, the map $\cyl(\phi_{X,i}):\cyl(X_i) \to \cyl(\liminj_i
X_i)$ is then the inclusion $S_i\p \{0,1\} \subset (\liminj_i S_i) \p
\{0,1\}$ on states and the inclusion $L_i\p \{0,1\} \subset (\liminj_i
L_i) \p \{0,1\}$ on actions. Thus, the map $\liminj_i \cyl(\phi_{X,i})
: \liminj_i \cyl(X_i) \to \cyl(\liminj_i X_i)$ induces a bijection on
states and on actions since the category of sets is cartesian-closed
(for the sequel, we will suppose that $\liminj_i \cyl(\phi_{X,i})$ is
the identity on states and on actions by abuse of
notation). Consequently, by
\cite[Proposition~4.4]{homotopyprecubical}, the map $\liminj_i
\cyl(\phi_{X,i}) : \liminj_i \cyl(X_i) \to \cyl(\liminj_i X_i)$ is
one-to-one on transitions. Let
$((\alpha,\epsilon_0),(u_1,\epsilon_1),\dots,(u_n,\epsilon_n),(\beta,\epsilon_{n+1}))$
be a transition of $\cyl(\liminj_i X_i)$. By definition of $\cyl$, the
tuple $(\alpha,u_1,\dots,u_n,\beta)$ is a transition of $\liminj_i
X_i$. Let $T_i$ be the image by the map $\phi_{X,i}:X_i\to \liminj_i
X_i$ of the set of transitions of $X_i$. Let $G_0(\bigcup T_i)=\bigcup
T_i$. Let us define $G_\lambda(\bigcup_{i} T_i)$ by induction on the
transfinite ordinal $\lambda\geq 0$ by $G_\lambda(\bigcup_{i} T_i)=
\bigcup_{\kappa<\lambda} G_\kappa(\bigcup_{i} T_i)$ for every limit
ordinal $\lambda$ and $G_{\lambda+1}(\bigcup_{i} T_i)$ is obtained
from $G_\lambda(\bigcup_{i} T_i)$ by adding to $G_\lambda(\bigcup_{i}
T_i)$ all tuples obtained by applying the patching axiom to tuples of
$G_\lambda(\bigcup_{i} T_i)$ in $\liminj_i X_i$. Hence we have the
inclusions $G_\lambda(\bigcup_{i} T_i) \subset
G_{\lambda+1}(\bigcup_{i} T_i)$ for all $\lambda\geq 0$. For
cardinality reason, there exists an ordinal $\lambda_0$ such that for
every $\lambda\geq \lambda_0$, there is the equality
$G_\lambda(\bigcup_{i} T_i) = G_{\lambda_0}(\bigcup_{i} T_i)$. The set
$G_{\lambda_0}(\bigcup_{i} T_i)$ is the set of transitions of
$\liminj_i X_i$ by \cite[Proposition~3.5]{hdts}. We are going to prove
by transfinite induction on $\lambda \geq 0$ the assertion
$\LL_\lambda$: \emph{if $(\alpha,u_1,\dots,u_n,\beta)\in
  G_\lambda(\bigcup_{i} T_i)$, then the tuple
  $((\alpha,\epsilon_0),(u_1,\epsilon_1),\dots,(u_n,\epsilon_n),(\beta,\epsilon_{n+1}))$
  is a transition of $\liminj_i \cyl(X_i)$ for any choice of
  $\epsilon_0,\dots,\epsilon_{n+1}\in \{0,1\}$}. Assume that $\lambda
= 0$. This implies that there exists a transition
$(\alpha^{i_0},u^{i_0}_1,\dots,u^{i_0}_n,\beta^{i_0})$ of some
$X_{i_0}$ such that
$\phi_{X,i_0}(\alpha^{i_0},u^{i_0}_1,\dots,u^{i_0}_n,\beta^{i_0})=(\alpha,u_1,\dots,u_n,\beta)$.
In particular, this means that $\phi_{X,i_0}(\alpha^{i_0})=\alpha$,
$\phi_{X,i_0}(\beta^{i_0})=\beta$ and for all $1\leq i \leq n$,
$\phi_{X,i_0}(u^{i_0}_i)=u_i$.  By definition of the functor $\cyl$,
we obtain
$\cyl(\phi_{X,i_0})(\alpha^{i_0},\epsilon_0)=(\alpha,\epsilon_0)$,
$\cyl(\phi_{X,i_0})(\beta^{i_0},\epsilon_{n+1})=(\beta,\epsilon_{n+1})$
and for all $1\leq i \leq n$,
$\cyl(\phi_{X,i_0})(u^{i_0}_i,\epsilon_i)=(u_i,\epsilon_i)$.  Since we
have $(\liminj_i \cyl(\phi_{X,i})) \phi_{\cyl X,i_0} =
\cyl(\phi_{X,i_0})$ by the universal property of the colimit, we
obtain $\phi_{\cyl
  X,i_0}(\alpha^{i_0},\epsilon_0)=(\alpha,\epsilon_0)$, $\phi_{\cyl
  X,i_0}(\beta^{i_0},\epsilon_{n+1})=(\beta,\epsilon_{n+1})$ and for
all $1\leq i \leq n$, $\phi_{\cyl
  X,i_0}(u^{i_0}_i,\epsilon_i)=(u_i,\epsilon_i)$. However, the tuple
$((\alpha^{i_0},\epsilon_0),(u^{i_0}_1,\epsilon_1),\dots,(u^{i_0}_n,\epsilon_n),(\beta^{i_0},\epsilon_{n+1}))$
is a transition of $\cyl(X_{i_0})$ by definition of the functor
$\cyl$. This implies that
\begin{multline*}
\phi_{\cyl X,i_0}((\alpha^{i_0},\epsilon_0),(u^{i_0}_1,\epsilon_1),\dots,(u^{i_0}_n,\epsilon_n),(\beta^{i_0},\epsilon_{n+1})) \\= ((\alpha,\epsilon_0),(u_1,\epsilon_1),\dots,(u_n,\epsilon_n),(\beta,\epsilon_{n+1}))
\end{multline*}
is a transition of $\liminj_i \cyl(X_i)$. We have proved $\LL_0$. Assume
$\LL_\kappa$ proved for all $\kappa<\lambda$ for some limit ordinal
$\lambda$. If $(\alpha,u_1,\dots,u_n,\beta)\in G_\lambda(\bigcup_{i}
T_i)$, then $(\alpha,u_1,\dots,u_n,\beta)\in G_\kappa(\bigcup_{i}
T_i)$ for some $\kappa<\lambda$, and therefore the tuple
\[((\alpha,\epsilon_0),(u_1,\epsilon_1),\dots,(u_n,\epsilon_n),(\beta,\epsilon_{n+1}))\]
is a transition of $\liminj_i \cyl(X_i)$ as well by induction
hypothesis. We have proved $\LL_\lambda$. Assume $\LL_\lambda$ proved
for $\lambda \geq 0$ and assume that $(\alpha,u_1,\dots,u_n,\beta)$
belongs to \[G_{\lambda+1}(\bigcup_{i} T_i) \backslash
G_{\lambda}(\bigcup_{i} T_i).\] Then there exist five tuples
\begin{align*}
&(\alpha',u'_1,\dots,u'_{n'},\beta')\\ &(\alpha',u'_1,\dots,u'_p,\nu'_1)\\&
(\nu'_1,u'_{p+1},\dots,u'_{n'},\beta')\\&(\alpha',u'_1,\dots,u'_{p+q},\nu'_2)\\&
(\nu'_2,u'_{p+q+1},\dots,u'_{n'},\beta')
\end{align*}
 of $G_{\lambda}(\bigcup_{i}
T_i)$ such that $(\nu'_1,u'_{p+1},\dots,u'_{p+q},\nu'_2) =
(\alpha,u_1,\dots,u_n,\beta)$. By induction hypothesis, the five
tuples
\begin{align*}
&((\alpha',0),(u'_1,\epsilon'_1),\dots,(u'_{n'},\epsilon_{n'}),(\beta',0))
\\&((\alpha',0),(u'_1,\epsilon'_1),\dots,(u'_p,\epsilon'_p),(\nu'_1,\epsilon_0))\\&
((\nu'_1,\epsilon_0),(u'_{p+1},\epsilon'_{p+1}),\dots,(u'_{n'},\epsilon'_{n'}),(\beta',0))\\
&((\alpha',0),(u'_1,\epsilon'_1),\dots,(u'_{p+q},\epsilon'_{p+q}),(\nu'_2,\epsilon_{n+1}))\\&
((\nu'_2,\epsilon_{n+1}),(u'_{p+q+1},\epsilon'_{p+q+1}),\dots,(u'_{n'},\epsilon'_{n'}),(\beta',0))
\end{align*}
are transitions of $\liminj_i \cyl(X_i)$ for any choice of
$\epsilon'_i \in \{0,1\}$. Therefore the tuple 
\[((\nu'_1,\epsilon_0),(u'_{p+1},\epsilon'_{p+1}),\dots,(u'_{p+q},\epsilon'_{p+q}),(\nu'_2,\epsilon_{n+1}))
\]
is a transition of $\liminj_i \cyl(X_i)$ by applying the patching
axiom in $\liminj_i \cyl(X_i)$. Let $\epsilon'_i = \epsilon_{i-p}$ for
$p + 1 \leq i \leq p + n$ and $\epsilon'_i=0$ otherwise. Since there
is the equality
\begin{multline*}
((\nu'_1,\epsilon_0),(u'_{p+1},\epsilon'_{p+1}),\dots,(u'_{p+q},\epsilon'_{p+q}),(\nu'_2,\epsilon_{n+1}))
\\=
((\alpha,\epsilon_0),(u_1,\epsilon_1),\dots,(u_n,\epsilon_n),(\beta,\epsilon_{n+1})),\end{multline*}
we deduce that $\LL_{\lambda+1}$ holds. The transfinite induction is
complete. We have proved that $\liminj_i \cyl(\phi_{X,i}): \liminj_i
\cyl(X_i) \to \cyl(\liminj_i X_i)$ is onto on transitions. The latter
map is bijective on states, bijective on actions and bijective on
transitions: it is an isomorphism of weak transition systems and the
proof is complete.  \epf

\bp \label{cocyl-def} Let $X=(S,\mu:L\to \Sigma,T)$ be a weak
transition system. There exists a well-defined weak transition system
$\cocyl(X)$ such that:
\begin{itemize}
\item The set of states is the set $S\p S$.
\item The set of actions is the set $L\p_\Sigma L$ and the labelling map 
is the canonical map $L\p_\Sigma L \to \Sigma$.
\item The transitions are the tuples
  $((\alpha^0,\alpha^1),(u_1^0,u_1^1),\dots,(u_n^0,u_n^1),(\beta^0,\beta^1))$
  such that for any $\epsilon_0,\dots,\epsilon_{n+1}\in \{0,1\}$, the
  tuple $(\alpha^{\epsilon_0},u_1^{\epsilon_1},\dots
  ,u_n^{\epsilon_n}, \beta^{\epsilon_{n+1}})$ is a transition of $X$.
\end{itemize}
Let $f:X\to Y$ be a map of weak transition systems. There exists a map
of weak transition systems $\cocyl(f):\cocyl(X) \to \cocyl(Y)$ defined
on states by the mapping $(\alpha^0,\alpha^1)\mapsto
(f(\alpha^0),f(\alpha^1))$ and on actions by the mapping
$(u^0,u^1)\mapsto (f(u^0),f(u^1))$.  
\ep

\bpf Let
$((\alpha^0,\alpha^1),(u_1^0,u_1^1),\dots,(u_n^0,u_n^1),(\beta^0,\beta^1))$
be a transition of $\cocyl(X)$. Let $\sigma$ be a permutation of
$\{1,\dots,n\}$ with $n\geq 2$. Then for any
$\epsilon_0,\dots,\epsilon_{n+1}\in \{0,1\}$, the tuple
$(\alpha^{\epsilon_0},u_{\sigma(1)}^{\epsilon_1},\dots
,u_{\sigma(n)}^{\epsilon_n}, \beta^{\epsilon_{n+1}})$ is a transition
of $X$ by the multiset axiom. Thus, the tuple
$((\alpha^0,\alpha^1),(u_{\sigma(1)}^0,u_{\sigma(1)}^1),\dots,(u_{\sigma(n)}^0,u_{\sigma(n)}^1),(\beta^0,\beta^1))$
is a transition of $\cocyl(X)$. Let $n\geq 3$. Let $p,q\geq 1$ with
$p+q<n$. Suppose that the five tuples
\begin{align*}
&((\alpha^0,\alpha^1),(u_1^0,u_1^1),\dots,(u_n^0,u_n^1),(\beta^0,\beta^1))\\
&((\alpha^0,\alpha^1),(u_1^0,u_1^1),\dots,(u_p^0,u_p^1),(\nu_1^0,\nu_1^1))\\
&((\nu_1^0,\nu_1^1),(u_{p+1}^0,u_{p+1}^1),\dots,(u_n^0,u_n^1),(\beta^0,\beta^1))\\
&((\alpha^0,\alpha^1),(u_1^0,u_1^1),\dots,(u_{p+q}^0,u_{p+q}^1),(\nu_2^0,\nu_2^1))\\
&((\nu_2^0,\nu_2^1),(u_{p+q+1}^0,u_{p+q+1}^1),\dots,(u_n^0,u_n^1),(\beta^0,\beta^1))
\end{align*}
are transitions of $\cocyl(X)$.  Then for any
$\epsilon_0,\dots,\epsilon_{n+1}\in \{0,1\}$, the tuple
\[(\nu_1^{\epsilon_0}, u_{p+1}^{\epsilon_{p+1}}, \dots,
u_{p+q}^{\epsilon_{p+q}}, \nu_2^{\epsilon_{n+1}})\] is a transition of
$X$ by the patching axiom. Thus, the tuple
\[((\nu_1^0,\nu_1^1),(u_{p+1}^0,u_{p+1}^1),\dots,(u_{p+q}^0,u_{p+q}^1),(\nu_2^0,\nu_2^1))\]
is a transition of $\cocyl(X)$. Hence $\cocyl(X)$ is well-defined as a
weak transition system. Let $f:X\to Y$ be a map of weak transition
systems. For any state $(\alpha^0,\alpha^1)$ of $\cocyl(X)$, the pair
$(f(\alpha^0),f(\alpha^1))$ is a state of $\cocyl(Y)$ by definition of
the functor $\cocyl$.  For any state $(u^0,u^1)$ of $\cocyl(X)$, we
have $\mu(u^0)=\mu(u^1)$ by definition of the functor $\cocyl$. We
deduce that $\mu(f(u^0))=\mu(u^0)=\mu(u^1)=\mu(f(u^1))$. Hence the
pair $(f(u^0),f(u^1))$ is an action of $\cocyl(Y)$ by definition of
the functor $\cocyl$.  Let
\[((\alpha^0,\alpha^1),(u_1^0,u_1^1),\dots,(u_n^0,u_n^1),(\beta^0,\beta^1))\]
be a transition of $\cocyl(X)$. By definition of the functor $\cocyl$,
the tuple
\[(\alpha^{\epsilon_0},u_1^{\epsilon_1},\dots,u_n^{\epsilon_n},\beta^{\epsilon_{n+1}})\]
is a transition of $X$ for any choice of
$\epsilon_0,\dots,\epsilon_{n+1}\in \{0,1\}$. Consequently, the tuple
\[(f(\alpha^{\epsilon_0}),f(u_1^{\epsilon_1}),\dots,f(u_n^{\epsilon_n}),f(\beta^{\epsilon_{n+1}}))\]
is a transition of $Y$ for any choice of
$\epsilon_0,\dots,\epsilon_{n+1}\in \{0,1\}$. By definition of the
functor $\cocyl$, we deduce that the tuple
\[((f(\alpha^0),f(\alpha^1)),(f(u_1^0),f(u_1^1)),\dots,(f(u_n^0),f(u_n^1)),(f(\beta^0),f(\beta^1)))\]
is a transition of $\cocyl(Y)$. We have proved the last part of the
statement.  \epf

We obtain a well-defined functor $\cocyl:\wts \to \wts$.  For
$\epsilon \in \{0,1\}$, there exists a unique map of weak transition
systems $\pi_X^\epsilon:\cocyl(X) \to X$ induced by the mappings
$(\alpha^0,\alpha^1)\mapsto \alpha^\epsilon$ on states and
$(u^0,u^1)\mapsto u^\epsilon$ on actions. Let
$\pi_X=(\pi_X^0,\pi_X^1)$. This defines a natural transformation $\pi
: \cocyl\Rightarrow \id\p \id$.

Since $\wts$ is locally presentable, and since the functor $\cyl:\wts
\to \wts$ is colimit-preserving by Proposition~\ref{cyl-colim}, we can
deduce that it is a left adjoint by applying the opposite of the
Special Adjoint Functor Theorem. The right adjoint is calculated in
the following proposition.

\bp There is a natural bijection of sets \[\Phi:\wts(\cyl(X),X')
\stackrel{\iso} \longrightarrow \wts(X,\cocyl(X'))\] for any weak
transition systems $X$ and $X'$.\ep

\bpf The proof is in seven parts.

\emph{1) Construction of $\Phi$.} Let $X=(S,\mu:L\to \Sigma,T)$ and
$X'=(S',\mu:L'\to \Sigma,T')$ be two weak transition systems.  let
$f\in \wts(\cyl(X),X')$. Let $g^0:S\to S'\p S'$ be the set map defined
by $g^0(\alpha)=(f^0(\alpha,0),f^0(\alpha,1))$. Let $\widetilde{g}:L
\to L'\p_\Sigma L'$ be the set map defined by
$\widetilde{g}(u)=(\widetilde{f}(u,0),\widetilde{f}(u,1))$.  Let
$(\alpha,u_1,\dots,u_n,\beta)$ be a transition of $X$. Then for any
$\epsilon_0,\dots,\epsilon_{n+1}\in \{0,1\}$, the tuple
$((\alpha,\epsilon_0),(u_1,\epsilon_1),\dots,(u_n,\epsilon_n),(\beta,\epsilon_{n+1}))$
is a transition of $\cyl(X)$ by definition of the functor
$\cyl$. Thus, the tuple
\[(f(\alpha,\epsilon_0),f(u_1,\epsilon_1),\dots,f(u_n,\epsilon_n),f(\beta,\epsilon_{n+1}))\]
is a transition of $X'$ since $f$ is a map of weak transition
systems. We deduce that the tuple
\[((f(\alpha,0),f(\alpha,1)),(f(u_1,0),f(u_1,1)),\dots,(f(u_n,0),f(u_n,1)),(f(\beta,0),f(\beta,1))\]
is a transition of $\cocyl(X')$ by definition of $\cocyl$. We have obtained a natural set map 
\[g=\Phi(f) : \wts(\cyl(X),X')
\longrightarrow \wts(X,\cocyl(X')).\] 

\emph{2) The case $X=\varnothing$.} There is the equality 
$\cyl(\varnothing)=\varnothing$.  We obtain the bijection
$\wts(\cyl(\varnothing),X')\iso \wts(\varnothing,\cocyl(X'))$.  We
have proved that $\Phi$ induces a bijection for $X=\varnothing$.

\emph{3) The case $X=\{0\}$.} There is the equality
\[\wts(\cyl(\{0\}),X') \iso \wts(\{(0,0),(0,1)\},X')\] by definition of
$\cyl$. There is the equality \[\wts(\{(0,0),(0,1)\},X')\iso
\wts(\{(0,0)\} \sqcup \{(0,1)\},X')\] by
\cite[Proposition~5.6]{cubicalhdts}. Hence we obtain the bijection
\[\wts(\cyl(\{0\}),X')\iso \wts(\{(0,0)\},X') \p \wts(\{(0,1)\},X').\]
The right-hand term is equal to $S'\p S'$, which is precisely
$\wts(\{0\},\cocyl(X'))$ by definition of $\cocyl$. We have proved
that $\Phi$ induces a bijection for $X=\{0\}$.

\emph{4) The case $X=\underline{x}$ for $x\in \Sigma$.}  There is the
equality $\cyl(\underline{x}) = \underline{x} \sqcup
\underline{x}$. Therefore we obtain the bijections
$\wts(\cyl(\underline{x}),X') \iso \wts(\underline{x} \sqcup
\underline{x},X')\iso \wts(\underline{x},X') \p
\wts(\underline{x},X')$.  The set $\wts(\cyl(\underline{x}),X')$ is
then equal to $\mu^{-1}(x)\p \mu^{-1}(x)$. And the set
$\wts(\underline{x},\cocyl(X'))$ is the set of actions of $\cocyl(X')$
labelled by $x$, i.e. $\mu^{-1}(x)\p \mu^{-1}(x)$.  We have proved
that $\Phi$ induces a bijection for $X=\underline{x}$ for all $x\in
\Sigma$.

\emph{5) The case $X=C_n^{ext}[x_1,\dots,x_n]$.} The
set of transitions of $\cyl(C_n^{ext}[x_1,\dots,x_n])$ is the set of
tuples
$((0_n,\epsilon_0),((x_{\sigma(1)},\sigma(1)),\epsilon_1),\dots,((x_{\sigma(n)},\sigma(n)),\epsilon_n),(1_n,\epsilon_{n+1}))$
for $\epsilon_0,\dots,\epsilon_{n+1} \in \{0,1\}$ and all permutation
$\sigma$ of $\{1,\dots,n\}$. A map \[f:\cyl(C_n^{ext}[x_1,\dots,x_n]) \longrightarrow X'\] is then determined by the
choice of four states $f(0_n,0),f(0_n,1),f(1_n,0),f(1_n,1)$ of $X'$
and for every $1\leq i\leq n$ by the choice of two actions
$f((x_i,i),0)$ and $f((x_i,i),1)$ of $X'$ such that the tuples \[
(f(0_n,\epsilon_0),f((x_{\sigma(1)},\sigma(1)),\epsilon_1),\dots,f((x_{\sigma(n)},\sigma(n)),\epsilon_n),f(1_n,\epsilon_{n+1}))
\] are transitions of $X'$ for all $\epsilon_0,\dots,\epsilon_{n+1} \in
\{0,1\}$ and all permutation $\sigma$ of $\{1,\dots,n\}$. By
definition of the functor $\cocyl$, the latter assertion is equivalent
to saying that the tuple
\begin{multline*}
((f(0_n,0),f(0_n,1)),(f((x_1,1),0),f((x_1,1),1)),\dots,\\(f((x_n,n),0),f((x_n,n),1)),
(f(1_n,0),f(1_n,1)))
\end{multline*}
is a transition of $\cocyl(X')$. Choosing a map $f$ from
$\cyl(C_n^{ext}[x_1,\dots,x_n])$ to $X'$ is therefore equivalent to
choosing a map of $\wts(C_n^{ext}[x_1,\dots,x_n],\cocyl(X'))$. We have
proved that $\Phi$ induces a bijection for
$X=C_n^{ext}[x_1,\dots,x_n]$ for $n\geq 1$ and for all $x_1,\dots,x_n
\in \Sigma$.

\emph{6) The case $X = X_1 \sqcup X_2$.} If $\Phi$ induces the bijections of sets 
$\wts(\cyl(X_i),X') \iso \wts(X_i,\cocyl(X'))$  for $i=1,2$, then we obtain the sequence of bijections
{\small
\begin{align*} 
&\wts(\cyl(X),X') & \\
 & \iso \wts(\cyl(X_1 \sqcup X_2),X') & \hbox{ by definition of $X$}\\
&\iso \wts(\cyl(X_1) \sqcup \cyl(X_2),X') & \hbox{ by Proposition~\ref{cyl-colim}}\\
&\iso \wts(\cyl(X_1),X') \p \wts(\cyl(X_2),X') & \hbox{ since $\wts(-,X')$ is limit-preserving}\\
&\iso \wts(X_1,\cocyl(X')) \p \wts(X_2,\cocyl(X')) & \hbox{ by hypothesis}\\
&\iso \wts(X_1 \sqcup X_2,\cocyl(X'))&  \hbox{ since $\wts(-,\cocyl(X'))$ is limit-preserving}\\
&\iso \wts(X,\cocyl(X')) & \hbox{ by definition of $X$}.
\end{align*}
}
We have proved that $\Phi$ induces a bijection for $X = X_1 \sqcup
X_2$.

\emph{7) End of the proof.} The functor $X\mapsto \wts(\cyl(X),X')$
from the opposite of the category $\wts$ to the category of sets is
limit-preserving by Proposition~\ref{cyl-colim}.  The functor
$X\mapsto \wts(X,\cocyl(X'))$ from the opposite of the category $\wts$
to the category of sets is limit-preserving as well since the functor
$\wts(-,Z)$ is limit-preserving as well for any weak transition system
$Z$. The proof is complete by observing that the canonical map
$\varnothing \to X$ belongs to $\cell_\wts(\I)$ by
Proposition~\ref{small-mono}.  \epf

\begin{cor} The weak transition system $V$ is exponential. \end{cor}

\bp Let $f:X\to X'$ be a monomorphism of $\wts$. Then the maps $f\star
\gamma^0$, $f\star \gamma^1$ and $f\star \gamma$ are monomorphisms of
$\wts$. \ep

\bpf Let $X=(S,\mu:L\to \Sigma,T)$ and $X'=(S',\mu:L'\to \Sigma,T')$.
The map $f\star \gamma^\epsilon$ induces on states the set map
$S'\sqcup_{S\p \{\epsilon\}} (S\p \{0,1\}) \longrightarrow S'\p
\{0,1\}$ which is one-to-one since the map $S\to S'$ is
one-to-one. And it induces on actions the set map $L'\sqcup_{L\p
  \{\epsilon\}} (L\p \{0,1\}) \longrightarrow L'\p \{0,1\}$ which is
one-to-one since the map $L\to L'$ is one-to-one. So by
\cite[Proposition~3.1]{cubicalhdts}, the map $f\star
\gamma^\epsilon:X'\sqcup_X \cyl(X) \to \cyl(X')$
is a monomorphism of $\wts$. The map $f\star \gamma$ induces on states
the set map $(S'\sqcup S')\sqcup_{S\sqcup S} (S\p \{0,1\})
\longrightarrow S'\p \{0,1\}$ which is the identity of $S'\sqcup
S'$. And it induces on actions the identity of $L'\to L'$. So by
\cite[Proposition~3.1]{cubicalhdts}, the map $f\star \gamma:(X'\sqcup
X')\sqcup_{X\sqcup X} \cyl(X) \to \cyl(X')$ is a
monomorphism of $\wts$.  \epf

\begin{cor}
  The cylinder $\cyl:\wts \to \wts$ is cartesian with respect
  to the class of monomorphisms of weak transition systems.
\end{cor}

We have all the ingredients leading to an Olschok model structure (see
\cite[Definition~2.7]{leftdet} for the definition of an Olschok model
structure):

\bth \label{standard-wts} There exists a unique left determined model
category on $\wts$ such that the cofibrations are the
monomorphisms. This model structure is an Olschok model structure,
with the very good cylinder $\cyl$ above defined.  \eth

\bpf This a consequence of Olschok's theorems.
\epf

\section{Restricting the model structure of weak transition systems} \label{restriction-modelstructure}

We start this section by restricting the model structure to the full
subcategory of cubical transition systems.

By definition, a \emph{cubical transition system} satisfies all axioms of
weak transition systems and the following two additional axioms (with
the notations of Definition~\ref{def_HDTS}):
\begin{itemize}
\item (All actions are used) For every $u\in L$, there is a
  $1$-transition $(\alpha,u,\beta)$.
\item (Intermediate state axiom) For every $n\geq 2$, every $p$ with
  $1\leq p<n$ and every transition $(\alpha,u_1,\dots,u_n,\beta)$ of
  $X$, there exists a state $\nu$ such that both
  $(\alpha,u_1,\dots,u_p,\nu)$ and $(\nu,u_{p+1},\dots,u_n,\beta)$ are
  transitions.
\end{itemize}

By definition, a cubical transition system is \emph{regular} if it
satisfies the Unique intermediate state axiom, also called CSA2:
 \begin{itemize}
\item (Unique intermediate state axiom or CSA2)  For every
  $n\geq 2$, every $p$ with $1\leq p<n$ and every transition
  $(\alpha,u_1,\dots,u_n,\beta)$ of $X$, there exists a unique state
  $\nu$ such that both $(\alpha,u_1,\dots,u_p,\nu)$ and
  $(\nu,u_{p+1},\dots,u_n,\beta)$ are transitions.
\end{itemize}

Here is an important example of regular transition systems: 
\begin{itemize}
\item For every $x\in \Sigma$, let us denote by $\dd{x}$ the cubical
  transition system with four states $\{1,2,3,4\}$, one action $x$ and
  two transitions $(1,x,2)$ and $(3,x,4)$. The cubical transition
  system $\dd{x}$ is called the \emph{double transition (labelled by
    $x$)} where $x\in \Sigma$.
\end{itemize}

\begin{nota} The full subcategory of $\wts$ of cubical transition
  systems is denoted by $\cts$.  The full subcategory of $\cts$ of
  regular transition systems is denoted by $\rts$.
\end{nota}

The category $\rts$ of regular transition systems is a reflective
subcategory of the category $\cts$ of cubical transition systems by
\cite[Proposition~4.4]{csts}. The reflection is denoted by
$\CSA_2:\cts \to \rts$.  The unit of the adjunction $\id \Rightarrow
\CSA_2$ forces CSA2 to be true by identifying the states provided by a
same application of the intermediate state axiom (see
\cite[Proposition~4.2]{csts}).

Let us introduce now the weak transition system corresponding to the
labelled $n$-cube.

\bp \label{cas_cube} \cite[Proposition~5.2]{hdts} Let $n\geq 0$ and
$x_1,\dots,x_n\in \Sigma$. Let $T_d\subset \{0,1\}^n \p
\{(x_1,1),\dots,(x_n,n)\}^d \p \{0,1\}^n$ (with $d\geq 1$) be the
subset of $(d+2)$-tuples
\[((\epsilon_1,\dots,\epsilon_n), (x_{i_1},i_1),\dots,(x_{i_d},i_d),
(\epsilon'_1,\dots,\epsilon'_n))\] such that
\begin{itemize}
\item $i_m = i_n$ implies $m = n$, i.e. there are no repetitions in the
  list \[(x_{i_1},i_1),\dots,(x_{i_d},i_d)\]
\item for all $i$, $\epsilon_i\leq \epsilon'_i$
\item $\epsilon_i\neq \epsilon'_i$ if and only if
  $i\in\{i_1,\dots,i_d\}$. 
\end{itemize}
Let $\mu : \{(x_1,1),\dots,(x_n,n)\} \rightarrow \Sigma$ be the set
map defined by $\mu(x_i,i) = x_i$. Then \[C_n[x_1,\dots,x_n] =
(\{0,1\}^n,\mu : \{(x_1,1),\dots,(x_n,n)\}\rightarrow
\Sigma,(T_d)_{d\geq 1})\] is a well-defined regular transition system
called the {\rm $n$-cube}. \ep

The $n$-cubes $C_n[x_1,\dots,x_n]$ for all $n\geq 0$ and all
$x_1,\dots,x_n\in \Sigma$ are regular by \cite[Proposition~5.2]{hdts}
and \cite[Proposition~4.6]{hdts}.  For $n = 0$, $C_0[]$, also denoted
by $C_0$, is nothing else but the one-state higher dimensional
transition system $(\{()\},\mu:\varnothing \rightarrow
\Sigma,\varnothing)$.

Since the tuple $(0,(x,0),0)$ is a transition of $V$ for all $x\in
\Sigma$, all actions are used. The intermediate state axiom is
satisfied since both the states $0$ or $1$ can always divide a
transition in two parts. Therefore the weak transition system $V$ is
cubical. Note that the cubical transition system $V$ is not regular.

\bth \label{standard-cts} There exists a unique left determined model
category on $\cts$ such that the cofibrations are the monomorphisms of
weak transition systems between cubical transition systems. This model
structure is an Olschok model structure, with the very good cylinder
$\cyl$ above defined.  \eth

\bpf The category $\cts$ is a full coreflective locally finitely
presentable subcategory of $\wts$ by
\cite[Corollary~3.15]{cubicalhdts}. The full subcategory of cubical
transition systems is a small injectivity class by
\cite[Theorem~3.6]{cubicalhdts}: more precisely being cubical is
equivalent to being injective with respect to the set of inclusions
$C_n[x_1,\dots,x_n]^{ext} \subset C_n[x_1,\dots,x_n]$ and
$\underline{x_1} \subset C_1[x_1]$ for all $n\geq 0$ and all
$x_1,\dots, x_n \in \Sigma$. Therefore, by
\cite[Proposition~4.3]{MR95j:18001}, it is closed under binary
products. Hence we obtain the inclusion $\cyl(\cts) \subset \cts$
since $V$ is cubical. Then \cite[Theorem~4.3]{leftdet} can be applied
because all maps $C_n[x_1,\dots,x_n]^{ext} \subset C_n[x_1,\dots,x_n]$
and $\underline{x_1} \subset C_1[x_1]$ for all $n\geq 0$ and all
$x_1,\dots, x_n \in \Sigma$ are cofibrations.  \epf

The right adjoint $\cocyl^{\cts}:\cts \to \cts$ of the restriction of
$\cyl$ to the full subcategory of cubical transition systems is the
composite map \[\cocyl^{\cts}:\cts \subset \wts \stackrel{\cocyl}
\longrightarrow \wts \longrightarrow \cts\] where the right-hand map
is the coreflection, obtained by taking the canonical colimit over
all cubes and all double transitions \cite[Theorem~3.11]{cubicalhdts}:
\[\cocyl^{\cts}(X) = \liminj_{{\begin{array}{c}f: C_n[x_1,\dots,x_n]
      \rightarrow \cocyl(X)\\\hbox{ or }f:\dd{x} \rightarrow \cocyl(X)\end{array}}}
\dom(f).\]
Therefore, we obtain:

\bp \label{inclusion} The counit map $\cocyl^{\cts}(X) \to \cocyl(X)$
is bijective on states and one-to-one on actions and transitions.  \ep

\bpf This is a consequence of the first part of
\cite[Theorem~3.11]{cubicalhdts}.  \epf

\begin{lem} \label{colim-preserv} The forgetful functor mapping a
  cubical transition system to its set of states is
  colimit-preserving.  The forgetful functor mapping a cubical
  transition system to its set of actions is colimit-preserving.
\end{lem}

\bpf Since the category of cubical transition systems is a
coreflective subcategory of the category of weak transition systems by
\cite[Corollary~3.15]{cubicalhdts}, this lemma is a consequence of
Lemma~\ref{colim-preserv0}.  \epf

Theorem~\ref{standard-cts} proves the existence of a set of generating
cofibrations for the model structure. It does not give any way to find
it.

\begin{lem} \label{R-cts-epic}
All maps of $\cell_\cts(\{R\})$ are epic.
\end{lem}

\bpf Let $f,g,h$ be three maps of $\cts$ with $f\in \cell_\cts(\{R\})$
such that $gf=hf$. Since $\cts$ is coreflective in $\wts$, we obtain
$f\in \cell_\wts(\{R\})$. Since $\cts$ is a full subcategory of
$\wts$, we obtain $gf=hf$ in $\wts$. By Lemma~\ref{R-wts-epic}, we
obtain $g=h$. \epf

\bth (Compare with \cite[Notation~4.5]{homotopyprecubical} and
\cite[Theorem~4.6]{homotopyprecubical}) \label{explicit} The set of
maps
\begin{multline*}
\I^{\cts} = \{C:\varnothing \to \{0\}\} \\ \cup \{\de C_n[x_1,\dots,x_n]
\to C_n[x_1,\dots,x_n]\mid \hbox{$n\geq 1$ and $x_1,\dots,x_n \in
  \Sigma$}\} \\ \cup \{C_1[x] \to \dd{x}\mid x\in \Sigma\}.
\end{multline*}
generates the class of cofibrations of the model structure of $\cts$.
\eth

\bpf By \cite[Theorem~4.6]{homotopyprecubical}, a cofibration between
cubical transition systems belongs to $\cell_{\cts}(\I^{\cts} \cup
\{R\})$ where $R:\{0,1\} \to \{0\}$ is the map identifying two states
since it is one-to-one on actions. Every map of $\I^{\cts}$ is
one-to-one on states. Therefore, there is the inclusion $\I^{\cts}
\subset \inj_{\cts}(\{R\})$.  Every map of $\cell_{\cts}(\{R\})$ is
epic by Lemma~\ref{R-cts-epic}.  By Theorem~\ref{move_bad}, every
cofibration $f$ then factors as a composite $f=f^+ f^-$ such that
$f^-\in \cell_{\cts}(\{R\})$ and $f^+\in \cell_{\cts}(\I^{\cts})$,
i.e. $R$ can be relocated at the beginning of the cellular
decomposition. Since the cofibration $f$ is also one-to-one on states
by definition of a cofibration, the map $f^-\in \cell_{\cts}(\{R\})$
is one-to-one on states as well. Therefore $f^-$ is trivial and there
is the equality $f=f^+$. We deduce that $f$ belongs to
$\cell_{\cts}(\I^{\cts})$. Conversely, by Lemma~\ref{colim-preserv},
every map of $\cell_{\cts}(\I^{\cts})$ is one-to-one on states and on
actions. Consequently, the class of cofibrations of $\cts$ is
$\cell_{\cts}(\I^{\cts})$. Since the underlying category $\cts$ is
locally presentable, every map of $\cof_{\cts}(\I^{\cts})$ is a
retract of a map of $\cell_{\cts}(\I^{\cts})$. Therefore every map of
$\cof_{\cts}(\I^{\cts})$ is one-to-one on states and on actions. We
obtain $\cof_{\cts}(\I^{\cts}) \subset \cell_{\cts}(\I^{\cts})$. Hence
we have obtained $\cof_{\cts}(\I^{\cts}) = \cell_{\cts}(\I^{\cts})$
and the proof is complete.  \epf

\bd Let $X$ be a weak transition system. A state $\alpha$ of $X$ is
{\rm internal} if there exists three transitions
\[
(\gamma,u_1,\dots,u_n,\delta),
(\gamma,u_1,\dots,u_p,\alpha),
(\alpha,u_{p+1},\dots,u_n,\delta)
\]
with $n\geq 2$ and $p\geq 1$. A state $\alpha$ is {\rm external} if it
is not internal. \ed

\begin{nota} The set of internal states of a weak transition system
  $X$ is denoted by $X^0_{int}$. The complement is denoted by
  $X^0_{ext}=X^0 \backslash X^0_{int}$. \end{nota}

An internal state cannot be initial or final. The converse is false.
Consider the amalgamated sum $C_1[x]*C_1[y]$ with $x,y\in \Sigma$ where
the final state of $C_1[x]$ is identified with the initial state of
$C_1[y]$: the intermediate state is not internal because
$C_1[x]*C_1[y]$ does not contain any $2$-transition.

\bp \label{path-regular} Let $X$ be a regular transition system. Then
the cubical transition system $\cocyl^{\cts}(X)$ is regular.  \ep

\bpf Let $(\gamma^-,\gamma^+)$ and $(\delta^-,\delta^+)$ be two states
of $\cocyl^{\cts}(X)$ such that the four tuples
\begin{align*}
&((\alpha^-,\alpha^+),(u_1^-,u_1^+),\dots,(u_p^-,u_p^+),(\gamma^-,\gamma^+))\\&((\gamma^-,\gamma^+),(u_{p+1}^-,u_{p+1}^+),\dots,(u_n^-,u_n^+),(\beta^-,\beta^+))\\&
((\alpha^-,\alpha^+),(u_1^-,u_1^+),\dots,(u_p^-,u_p^+),(\delta^-,\delta^+))\\&
((\delta^-,\delta^+),(u_{p+1}^-,u_{p+1}^+),\dots,(u_n^-,u_n^+),(\beta^-,\beta^+))
\end{align*}
are transitions of $\cocyl^{\cts}(X)$, and therefore of $\cocyl(X)$ by
Proposition~\ref{inclusion}. By definition of $\cocyl(X)$, the tuples
\begin{align*}
&(\alpha^\pm,u_1^\pm,\dots,u_p^\pm,\gamma^\pm)
\\&(\gamma^\pm,u_{p+1}^\pm,\dots,u_n^\pm,\beta^\pm)\\&
(\alpha^\pm,u_1^\pm,\dots,u_p^\pm,\delta^\pm)\\&(\delta^\pm,u_{p+1}^\pm,\dots,u_n^\pm,\beta^\pm)
\end{align*}
are transitions of $X$. Since $X$ is regular, we obtain
$\gamma^-=\gamma^+=\delta^-=\delta^+$. In particular, this implies
that $(\gamma^-,\gamma^+)=(\delta^-,\delta^+)$. Hence the cubical
transition system $\cocyl^{\cts}(X)$ is regular as well.  \epf

\begin{lem} Let $X = (S,\mu:L\rightarrow \Sigma,T)$ be a weak
  transition system. Let $S'\subset S$. Let $T\restr S'$ be
  the subset of tuples of $T$ such that both the initial and the final
  states belong to $S'$. Then the triple $(S',\mu:L\rightarrow
  \Sigma,T\restr S')$ yields a well-defined weak transition
  system denoted by $X\restr S'$.
\end{lem}

\bpf For every permutation $\sigma$ of $\{1,\dots,n\}$ with $n\geq 2$,
if the tuple $(\alpha,u_1,\dots,u_n,\beta)$ is a transition such that
$\alpha,\beta \in S'$, then $(\alpha,u_{\sigma(1)}, \dots,
u_{\sigma(n)}, \beta) \in T\restr S'$. Therefore the set of tuples
$T\restr S'$ satisfies the multiset axiom.  For every $(n+2)$-tuple
$(\alpha,u_1,\dots,u_n,\beta)$ with $n\geq 3$, for every $p,q\geq 1$
with $p+q<n$, if the five tuples $(\alpha,u_1, \dots, u_n, \beta)$,
$(\alpha,u_1, \dots, u_p, \nu_1)$, $(\nu_1, u_{p+1}, \dots, u_n,
\beta)$, $(\alpha, u_1, \dots, u_{p+q}, \nu_2)$ and $(\nu_2,
u_{p+q+1}, \dots, u_n, \beta)$ are transitions of $T\restr S'$, then
$\nu_1,\nu_2 \in S'$. Therefore the transition $(\nu_1, u_{p+1},
\dots, u_{p+q}, \nu_2)$ belongs to $T\restr S'$. Thus, the set of
tuples $T\restr S'$ satisfies the patching axiom.  \epf

\begin{lem} \label{identifying-states-in-cyl} Let $X$ be a weak
  transition system. Let $Z \subset X^0$ be a subset of the set $X^0$
  of states of $X$. Consider a map $f:\cyl(X)\to Y$ of $\wts$ such
  that $f$ belongs to $\cell_\wts(\{R\})$ where $R:\{0,1\} \to \{0\}$
  is the set map identifying two states. Suppose that every cell of
  $f$ is of the form
\[
\xymatrix
{
\{0,1\} \fR{m_\alpha:\epsilon \mapsto (\alpha,\epsilon)} \fD{R} && \bullet \fD{} \\
&& \\
\{0\} \fR{} && \cocartesien \bullet,
}
\]
i.e. the identifications of states made by $f$ are of the form
$(\alpha,0)=(\alpha,1)$ for $\alpha\in Z$. Then $f$ is onto on maps,
bijective on actions, onto on transitions and split epic. There is an
inclusion $Y\subset \cyl(X)$ which is a section. Moreover, if $X$ is
cubical, then $Y$ is cubical as well.
\end{lem}

In general, identifications of states may generate new transitions by
the patching axiom.  The point is that it is not the case for this
particular situation.

\begin{nota} With the notations of
  Lemma~\ref{identifying-states-in-cyl} , let $Y=\cyl(X)//Z$.
\end{nota}

\bpf[Proof of Lemma~\ref{identifying-states-in-cyl}] The map
$R:\{0,1\} \to \{0\}$ is onto on states and bijective on
actions. Therefore the map $f$ is onto on states and bijective on
actions by Lemma~\ref{colim-preserv0}. Consider the cocone of
$\set^{\{s\}\cup \Sigma}$ consisting of the unique map
\[\omega(\cyl(X)) \longrightarrow (Z\p \{0\} \sqcup ((X^0\backslash
Z)\p \{0,1\}),L_X\p \{0,1\})\] where $L_X$ is the set of actions of
$X$.  The $\omega$-final lift gives rise to the map of weak transition
systems $f:\cyl(X) \to \cyl(X)//Z$. The final structure contains all
tuples of the form
$((\alpha,\epsilon_0),(u_1,\epsilon_1),\dots,(u_n,\epsilon_n),(\beta,\epsilon_{n+1}))$
such that $(\alpha,u_1,\dots,u_n,\beta)$ is a transition of $X$.  The
weak transition system $X\restr (Z\p \{0\} \sqcup ((X^0\backslash Z)\p
\{0,1\}))$ has exactly this set of transitions. Hence this set of
transitions is the final structure and $X\restr (Z\p \{0\} \sqcup
((X^0\backslash Z)\p \{0,1\})) = \cyl(X)//Z$. The identity on states
and the identity on actions induce a section of $f$, actually the
inclusion $\cyl(X)//Z\subset \cyl(X)$. Suppose moreover that $X$ is
cubical. Then the weak transition $\cyl(X)$ is cubical by
Theorem~\ref{standard-cts}.  Since the cocone above induces the
identity on actions, $\cyl(X)//Z$ is then cubical by
\cite[Theorem~3.3]{csts}.  \epf

\begin{lem} \label{calcul-new-cyl} Let $X$ be a regular transition system.
Then there is the natural isomorphism \[\CSA_2(\cyl(X)) \iso \cyl(X)//X^0_{int}.\]
\end{lem}

\bpf By Lemma~\ref{identifying-states-in-cyl}, the weak transition
system $\cyl(X)//X^0_{int}$ is cubical.  Let $(\mu_1,\zeta_1)$ and
$(\mu_2,\zeta_2)$ be two states of $\cyl(X)//X^0_{int}$ such that
there exists a transition
\[((\alpha,\epsilon_0),(u_1,\epsilon_1),\dots,(u_n,\epsilon_n),(\beta,\epsilon_{n+1}))\]
of $\cyl(X)//X^0_{int}$ such that the four tuples
\begin{multline*}
((\alpha,\epsilon_0),(u_1,\epsilon_1),\dots,(u_p,\epsilon_p),(\mu_1,\zeta_1))\\((\alpha,\epsilon_0),(u_1,\epsilon_1),\dots,(u_p,\epsilon_p),(\mu_2,\zeta_2))\\
((\mu_1,\zeta_1),(u_{p+1},\epsilon_{p+1}),\dots,(u_n,\epsilon_n),(\beta,\epsilon_{n+1}))\\
((\mu_2,\zeta_2),(u_{p+1},\epsilon_{p+1}),\dots,(u_n,\epsilon_n),(\beta,\epsilon_{n+1}))
\end{multline*}
are transitions of $\cyl(X)//X^0_{int}$ as well. Then the five tuples 
\begin{multline*}
(\alpha,u_1,\dots,u_n,\beta),
(\alpha,u_1,\dots,u_p,\mu_1),(\alpha,u_1,\dots,u_p,\mu_2)\\
(\mu_1,u_{p+1},\dots,u_n,\beta),
(\mu_2,u_{p+1},\dots,u_n,\beta)
\end{multline*}
are transitions of $X$ by definition of $\cyl(X)$.  Since $X$ is
regular, there is the equality $\mu_1=\mu_2$. Moreover, the state
$\mu_1=\mu_2$ belongs to $X^0_{int}$. Therefore $\zeta_1=\zeta_2=0$
and $(\mu_1,\zeta_1)=(\mu_2,\zeta_2)$. We deduce that
$\cyl(X)//X^0_{int}$ is regular. Consider a map $\cyl(X)\to Z$ with
$Z$ regular. The map $\omega(\cyl(X)) \to \omega(Z)$ makes the
identifications $(u,0)=(u,1)$ for all $u\in X^0_{int}$ by
CSA2. Therefore it factors uniquely as a composite \[\omega(\cyl(X))
\to (X^0_{int}\p\{0\} \sqcup ( X^0_{ext}\p \{0,1\}),L_X\p \{0,1\}) \to
\omega(Z).\] Hence the map $\cyl(X) \to Z$ factors uniquely as a
composite \[\cyl(X) \longrightarrow \cyl(X)//X^0_{int} \longrightarrow
Z.\] However, $\cyl(X)//X^0_{int}$ is regular. Hence we obtain the
isomorphism
\[\cyl(X)//X^0_{int}\iso \CSA_2(\cyl(X)).\]  \epf

\bp \label{section1} Let $X$ be a regular transition system. Then the
map \[\eta_{\cyl(X)}: \cyl(X) \to \CSA_2(\cyl(X))\] has a section
$s_{X}: \CSA_2(\cyl(X)) \to \cyl(X)$. \ep

\bpf By Lemma~\ref{calcul-new-cyl}, the inclusion $s_{X}:
\CSA_2(\cyl(X)) \subset \cyl(X)$ is a section of the natural map
$\eta_{\cyl(X)}: \cyl(X) \to \CSA_2(\cyl(X))$.  \epf

We can now prove the 

\bth \label{standard-rts} There exists a unique left determined model
category on $\rts$ such that the set of generating cofibrations is
$\CSA_2(\I^{\cts})=\I^{\cts}$. This model structure is an Olschok
model structure with the very good cylinder $\CSA_2 \cyl$.  \eth

\bpf Thanks to Proposition~\ref{path-regular} and
Proposition~\ref{section1}, the theorem is a consequence of
\cite[Theorem~3.1]{leftdet}.  \epf

For any regular transition system $X$, the map $\gamma_X:X\sqcup X \to
\CSA_2(\cyl(X))$ is a cofibration of the left determined model
structure of $\rts$. If $\alpha$ is an internal state of $X$, then the
two states $(\alpha,0)$ and $(\alpha,1)$ of $X \sqcup X$ are
identified in $\CSA_2(\cyl(X))$ by $\gamma_X$ by
Lemma~\ref{calcul-new-cyl}.  Consequently, as soon as $X$ contains
internal states, the cofibration $\gamma_X:X\sqcup X \to
\CSA_2(\cyl(X))$ is not one-to-one on states.

\section{The fibrant replacement functor destroys the causal structure} \label{boum}

We are going to prove in this section that the model structure of weak
transition systems as well as all its restrictions interact extremely
badly with the causal structure of the higher dimensional transition
systems. More precisely, the fibrant replacement functor destroys the
causal structure.

For the three model structures (on $\wts$, on $\cts$ and $\rts$), we
start from a weak transition system $X$ containing at least one
transition. We then consider the fibrant replacement $X^{fib}$ of $X$
in $\wts$ (in $\cts$ or in $\rts$ resp.) by factoring the canonical
map $X\to \mathbf{1}$ as a composite
\[\xymatrix{X \fR{\simeq}_-{i_X} \ar@{^(->}[rr] &&  X^{fib} \ar@{->>}[rr] && \mathbf{1}} \]
in $\wts$ (in $\cts$ or in $\rts$ resp.).  The canonical map
$X^{fib}\to \mathbf{1}$ is a fibration: therefore it satisfied the RLP
with respect to any trivial cofibration, in particular with respect to
any cofibration of the form $f\star \gamma^0$ where $f:A\to B$ is a
cofibration. By adjunction, the lift $\ell$ in any commutative diagram
of solid arrows of the form
\[
\xymatrix
{
A \fD{f} \fR{\phi} && P(X^{fib}) \fD{\pi^0}\\
&& \\
B \fR{\psi} \ar@{-->}[rruu]^-{\ell}&& X^{fib}
}
\]
exists, where $P(X^{fib})$ is the path space of $X^{fib}$ in $\wts$
(in $\cts$ or in $\rts$ resp.). Since $X$ contains at least one
transition, the image by $i_X$ gives rise to a transition
$(\alpha,u_1,\dots,u_n,\beta)$ of $X^{fib}$. Let us now treat first
the case of $\wts$, and then together the case of $\cts$ and
$\rts$. The key point in what follows is that, if $S$ denotes the set
of states of $X^{fib}$, then the cartesian product $S \p S$ is the set
of states of $P(X^{fib})$: for $\wts$, this is due to
Proposition~\ref{cocyl-def}, for $\cts$, this is a consequence of
Proposition~\ref{inclusion}, and finally for $\rts$, this is a
consequence of Proposition~\ref{path-regular}. The crucial fact is
that the coordinates in a cartesian product are independent from each
other.

\bth \label{path-wts} With the notations above, for any pair of states
$(\gamma,\delta)$ of the fibrant replacement $X^{fib}$ of $X$ in
$\wts$, the tuple $(\gamma,u_1,\dots,u_n,\delta)$ is a transition of
$X^{fib}$.  \eth

\bpf The transition $(\alpha,u_1,\dots,u_n,\beta)$ of $X^{fib}$ gives
rise to a map \[\psi:C_n^{ext}[\mu(u_1),\dots,\mu(u_n)]
\longrightarrow X^{fib}.\] We then consider the diagram above with the
cofibration \[f : \{0_n,1_n\} \subset
C_n^{ext}[\mu(u_1),\dots,\mu(u_n)]\] and with $\phi(0_n) =
(\alpha,\gamma)$ and $\phi(1_n) = (\beta,\delta)$. The existence of
the lift $\ell$ yields a transition
$((\alpha,\gamma),\ell(\mu(u_1),1),\dots,\ell(\mu(u_n),n),(\beta,\delta))$
of $P(X^{fib})$ with $\ell(\mu(u_i),i)=(u_i,u'_i)$ for some $u'_i$ for
$1\leq i \leq n$. By Proposition~\ref{cocyl-def}, we deduce that the
tuple $(\gamma,u_1,\dots,u_n,\delta)$ is a transition of $X^{fib}$.
\epf

Since the path functor in the category of cubical transition systems
is a subobject of the path functor in the category of weak transition
systems by Proposition~\ref{inclusion}, and since the path space in
$\cts$ of a regular transition system is regular by
Proposition~\ref{path-regular}, we can conclude in the same way:

\bth \label{path-cts} With the notations above. For any pair of states
$(\gamma,\delta)$ of the fibrant replacement $X^{fib}$ of $X$ in
$\cts$ (in $\rts$ resp.), the tuple $(\gamma,u_1,\dots,u_n,\delta)$ is
a transition of $X^{fib}$.  \eth

\bpf[Sketch of proof] Since we work now in $\cts$ (in $\rts$ resp.),
we must start from the cofibration $f : \{0_n,1_n\} \subset
C_n[\mu(u_1),\dots,\mu(u_n)]$ (the weak transition
$C_n^{ext}[\mu(u_1),\dots,\mu(u_n)]$ is not cubical nor regular since
it does not satisfy the intermediate state axiom).  The transition
$(\alpha,u_1,\dots,u_n,\beta)$ of $X^{fib}$ gives rise to a map
\[\psi:C_n^{ext}[\mu(u_1),\dots,\mu(u_n)] \longrightarrow X^{fib}.\] By
\cite[Theorem~3.6]{cubicalhdts}, the map $\psi$ factors as a composite
 \[\psi:C_n^{ext}[\mu(u_1),\dots,\mu(u_n)]\longrightarrow C_n[\mu(u_1),\dots,\mu(u_n)]
 \stackrel{\overline{\psi}}\longrightarrow X^{fib}.\] The rest of the
 proof is mutatis mutandis the proof of Theorem~\ref{path-wts}.
\epf

\section{The homotopy theory of star-shaped transition systems}

We need first to introduce some definitions and notations. In this
section, $\K$ is one of the three model categories $\wts$, $\cts$ or
$\rts$ equipped with the left determined model structure constructed
in this paper. Consider the one-state weak cubical regular transition
system ${\{\iota\}}$.  The forgetful functor
$\omega^{\{\iota\}}:{\{\iota\}}\ddownarrow \K \to \K$ defined on
objects by $\omega^{\{\iota\}}({\{\iota\}}\to X) = X$ and on maps by
$\omega^{\{\iota\}}({\{\iota\}}\to f)=f$ is a right adjoint. The left
adjoint $\rho^{\{\iota\}}:\K \rightarrow {\{\iota\}}\ddownarrow \K$ is
defined on objects by $\rho^{\{\iota\}}(X) = ({\{\iota\}}\to
{\{\iota\}} \sqcup X)$ and on morphisms by $\rho^{\{\iota\}}(f) =
\id_{\{\iota\}} \sqcup f$.  The locally presentable category
${\{\iota\}}\ddownarrow \K$ is equipped with the model structure
described in \cite[Theorem~2.7]{undercat}: a map $f$ is a cofibration
(fibration, weak equivalence resp.) of $i\ddownarrow \K$ if and only
if $\omega^{\{\iota\}}(f)$ is a cofibration (fibration, weak
equivalence resp.)  of $\K$. For the sequel, it is important to keep
in mind that the forgetful functor $\omega^{\{\iota\}}: {\{\iota\}}
\ddownarrow \K \to \K$ preserves colimits of connected diagrams, in
particular pushouts and transfinite compositions.

\bth \label{comma-cat} Let $\K$ be $\wts$ or $\cts$ or $\rts$. Then the
model category ${\{\iota\}}\ddownarrow \K$ is an Olschok model
category and is left determined. \eth

\bpf The map $\gamma_{\{\iota\}}:{\{\iota\}}\sqcup {\{\iota\}} \to \cyl({\{\iota\}})$
is an isomorphism by Proposition~\ref{cyl-construction}. Therefore it
is epic.  We can then apply \cite[Theorem~5.8]{leftdet} to obtain an
Olschok model structure.  Let $\underline{C}:\K\to \K$ be the cylinder
functor. By \cite[Lemma~5.7]{leftdet}, there is the pushout diagram of
$\K$:
\[
\xymatrix
{
{\{\iota\}}\sqcup {\{\iota\}} \fR{} \fD{} && {\{\iota\}} \ar@{->}[dd]^-{\underline{C}_{\{\iota\}}(\iota \to X)} \\
&& \\
\underline{C}(X) \fR{} && \cocartesien \omega^{\{\iota\}}(\underline{C}_{\{\iota\}}(\iota \to X))
}
\]
where $\underline{C}_{\{\iota\}}:\{\iota\} \ddownarrow \K \to
{\{\iota\}} \ddownarrow \K$ is the cylinder functor of the comma
category ${\{\iota\}} \ddownarrow \K$. The map \[\underline{C}(X) \to
\omega^{\{\iota\}}(\underline{C}_{\{\iota\}}(\{\iota\} \to X))\]
consists of the identification $(\iota,0)=(\iota,1)$. In $\wts$, and
in $\cts$ which is coreflective in $\wts$, the latter map is the map
\[\cyl(X) \to \cyl(X) // \{\iota\}\] which has a section by
Lemma~\ref{identifying-states-in-cyl}. Colimits in $\rts$ are
calculated first by taking the colimit in $\cts$ and then by applying
the reflection $\CSA_2:\cts\to \rts$. Let $X\in \rts$. The map
\[\underline{C}(X) \to
\omega^{\{\iota\}}(\underline{C}_{\{\iota\}}(\{\iota\} \to X))\]
consisting of the identification $(\iota,0)=(\iota,1)$ in $\rts$ 
is then equal to the composite map 
\begin{multline*}
\CSA_2(\cyl(X))\iso \cyl(X)//X^0_{int} \to \cyl(X) // (X^0_{int}\cup
\{\iota\}) \\ \to \CSA_2\left(\cyl(X) // (X^0_{int}\cup
\{\iota\})\right).
\end{multline*}
By Lemma~\ref{identifying-states-in-cyl}, the weak transition system
$\cyl(X) // (X^0_{int}\cup \{\iota\})$ is cubical and the map
$\cyl(X)//X^0_{int} \to \cyl(X) // (X^0_{int}\cup \{\iota\})$ has a
section: the inclusion \[\cyl(X) // (X^0_{int}\cup \{\iota\}) \subset
\cyl(X)//X^0_{int}.\] Therefore, there exists a map \[\cyl(X) //
(X^0_{int}\cup \{\iota\}) \to \CSA_2(\cyl(X))\] which is one-to-one on
states. Since $\CSA_2(\cyl(X))$ is regular, the cubical transition
system $\cyl(X) // (X^0_{int}\cup \{\iota\})$ is regular as well by
\cite[Proposition~4.1]{csts}. Hence there is an isomorphism 
\[\cyl(X) // (X^0_{int}\cup
\{\iota\}) \iso \CSA_2\left(\cyl(X) // (X^0_{int}\cup \{\iota\})\right).\]
We have proved that  the map
\[\underline{C}(X) \to
\omega^{\{\iota\}}(\underline{C}_{\{\iota\}}(\{\iota\} \to X))\]
consisting of the identification $(\iota,0)=(\iota,1)$ in $\rts$ is
the map \[\cyl(X)//X^0_{int} \to \cyl(X) // (X^0_{int}\cup
\{\iota\})\] which has a section by
Lemma~\ref{identifying-states-in-cyl}. 

Thanks to \cite[Corollary~5.9]{leftdet}, we deduce that the cylinder
$\underline{C}_{\{\iota\}}$ is very good and that the Olschok model
structure is left determined for the three cases $\K=\wts$, $\K=\cts$
and $\K=\rts$.  \epf

It is usual in computer science to work in the comma category ${\{\iota\}}
\ddownarrow \K$ where the image of the state $\iota$ represents the
initial state of the process which is modeled. It then makes sense to
restrict to the states which are reachable from this initial state by
a path of transitions. Hence we introduce the following definitions:

\bd Let $X$ be a weak transition system and let $\iota$ be a state of
$X$. A state $\alpha$ of $X$ is {\rm reachable from $\iota$} if it is
equal to $\iota$ or if there exists a finite sequence of transitions
$t_i$ of $X$ from $\alpha_i$ to $\alpha_{i+1}$ for $0\leq i \leq n$
with $n\geq 0$, $\alpha_0 = \iota$ and $\alpha_{n+1} = \alpha$.  \ed

\bd A {\rm star-shaped weak (cubical regular resp.) transition system}
is an object ${\{\iota\}} \rightarrow X$ of the comma category ${\{\iota\}}
\ddownarrow \K$ such that every state of the underlying weak
transition system $X$ is reachable from $\iota$. The full subcategory of
${\{\iota\}} \ddownarrow \K$ of star-shaped weak (cubical regular resp.)
transition systems is denoted by $\K_{\bullet}$. \ed

\bp \label{R-onto-comma} Let $\K$ be $\wts$ or $\cts$. Every map of
$\cell_{\{\iota\} \ddownarrow \K}(\{\rho^{\{\iota\}}(R)\})$ is onto on
states and the identity on actions. Moreover, every map of
$\cell_{\{\iota\} \ddownarrow \K}(\{\rho^{\{\iota\}}(R)\})$ is
epic. \ep

Proposition~\ref{R-onto-comma} also holds for $\K=\rts$ with a
slightly different proof.

\bpf By adjunction, if $f\in \cell_{\{\iota\} \ddownarrow
  \K}(\{\rho^{\{\iota\}}(R)\})$, then $\omega^{\{\iota\}}(f) \in
\cell_{\K}(R)$. Hence the first assertion of the proposition is a
consequence of Lemma~\ref{colim-preserv0} for $\K=\wts$ and of
Lemma~\ref{colim-preserv} for $\K=\cts$. Let $f,g,h$ be three maps of
$\{\iota\} \ddownarrow \K$ such that $f\in\cell_{\{\iota\} \ddownarrow
  \K}(\{\rho^{\{\iota\}}(R)\})$ and $gf=hf$.  Then we have by
functoriality $\omega^{\{\iota\}}(g)\omega^{\{\iota\}}(f) =
\omega^{\{\iota\}}(h)\omega^{\{\iota\}}(f)$. By Lemma~\ref{R-wts-epic}
if $\K=\wts$ and by Lemma~\ref{R-cts-epic} if $\K=\cts$, we obtain
$\omega^{\{\iota\}}(g) = \omega^{\{\iota\}}(h)$. Thus, there is the
equality $g=h$ and the proof is complete.
\epf

\bp \label{corefl_bullet} Let $\K$ be $\wts$ or $\cts$ or $\rts$. The
category $\K_\bullet$ is a coreflective full subcategory of
$\{\iota\}\ddownarrow \K$.  \ep

\bpf[Sketch of proof] The coreflection is described in
\cite[Proposition~6.5]{leftdet} for $\K=\wts$. For $\K=\cts$ or
$\K=\rts$, the coreflection $\{\iota\}\ddownarrow \K \to \K_\bullet$
removes all states which are not reachable from $\iota$, it removes
all transitions not starting from a reachable state and it removes all
actions which are not used by reachable transitions. This functor
clearly takes a cubical (regular resp.)  transition system to a
cubical (regular resp.) one since the state(s) dividing a transition
of the image is (are) reachable from $\iota$.  \epf

\bp \label{cone} Let $\K$ be $\wts$ or $\cts$ or $\rts$. The category
$\K_\bullet$ is a small-cone injectivity class of
$\{\iota\}\ddownarrow \K$ such that the top of the cone is
$\rho^{\{\iota\}}(\alpha)=\{\iota,\alpha\}$, such that the cone
contains only cofibrations and also the map $\{\iota,\alpha\}\to
\iota$.  \ep

\bpf For $\K=\wts$, this is \cite[Proposition~6.6]{leftdet}. For
$\K=\cts$ or $\K=\rts$, and since a cubical (regular resp.) transition
system satisfies the intermediate state axiom, a state is reachable
from $\iota$ if and only if it is reachable from $\iota$ by a path of
$1$-dimensional transitions. The cone consists of the map
$\{\iota,\alpha\}\to \iota$ and of the inclusions of
$\{\iota,\alpha\}$ into the cubical transition systems \[\iota
\stackrel{t_1} \longrightarrow \bullet \longrightarrow \dots
\longrightarrow \bullet \stackrel{t_n} \longrightarrow \alpha\] for
all $n\geq 1$ and all $1$-transitions $t_1,\dots,t_n$ with the
labelling map $\id_\Sigma$.  \epf

\begin{cor} Let $\K$ be $\wts$ or $\cts$ or $\rts$. The category
  $\K_\bullet$ is a small-cone injectivity class of
  $\{\iota\}\ddownarrow \K$ such that the cone contains only maps
  which are one-to-one on actions.
\end{cor}

\begin{cor} \label{locpre_bullet} Let $\K$ be $\wts$ or $\cts$ or
  $\rts$. The category $\K_\bullet$ is locally presentable.
\end{cor}

\bpf Since the category $\K_{\bullet}$ is a small cone-injectivity
class by Proposition~\ref{cone}, it is accessible by
\cite[Proposition~4.16]{MR95j:18001}. Therefore it is locally
presentable because it is a full coreflective subcategory of a
cocomplete category.  \epf

\bth \label{cofibrantly-gen} Let $\K$ be $\wts$ or $\cts$ or
$\rts$. The class of cofibrations of $\{\iota\}\ddownarrow \K$ between
objects of $\K_{\bullet}$ is cofibrantly generated.  \eth

\bpf By Proposition~\ref{cone}, there exists a set of cofibrations
$\{g:\{\iota,\alpha\} \to P\}$ of $\{\iota\}\ddownarrow \K$ such that
$\K_{\bullet}$ is the subcategory of $\{\iota\}\ddownarrow \K$ of
objects which are injective with respect to $\{g:\{\iota,\alpha\} \to
P\} \cup \{\{\iota,\alpha\}\to \{\iota\}\}$.  Consider a commutative
square of $\{\iota\}\ddownarrow \K$ of the form
\[
\xymatrix
{\{\iota\}\sqcup A \fR{\phi} \fD{\rho^{\{\iota\}}(f)} && X \fD{g} \\
&& \\
\{\iota\}\sqcup B \fR{\psi} && Y
}
\]
where $f$ is a generating cofibration of $\K$ and where the map
$g:X\to Y$ is a map between star-shaped objects. For every state
$\alpha$ of $A$, and since $X$ is star-shaped, the composite map
$\{\iota,\alpha\} \to \{\iota\}\sqcup A \to X$ factors as a composite
$\{\iota,\alpha\} \stackrel{g_\alpha}\to P_\alpha \to X$ with
$g_\alpha \in \{g:\{\iota,\alpha\} \to P\} \cup \{\{\iota,\alpha\}\to
\{\iota\}\}$. We obtain the commutative diagram of
$\{\iota\}\ddownarrow \K$
\[
\xymatrix
{
\{\iota\} \sqcup A^0 \fD{g_\alpha} \fR{} && \{\iota\} \sqcup A \fR{} \fD{} && X \fR{} && Y\\ 
&& && && \\
\bigsqcup\limits_{\alpha\in A^0} P_\alpha \fR{} && \cocartesien \widehat{A} \ar@{-->}[rruu]^-{\exists !} &&&&
}
\]
where $A^0$ is the set of states of $A$ and where the sum
$\bigsqcup_{\alpha\in A^0}$ is taken in the category
$\{\iota\}\ddownarrow \K$. The lift $\widehat{A}\to X$ exists and is
unique by the universal property of the pushout.  All generating
cofibrations of $\wts$ described in Notation~\ref{defI} and all
generating cofibrations of $\cts$ and $\rts$ described in
Theorem~\ref{explicit} and in Theorem~\ref{standard-rts} respectively
are one-to-one on states and on actions. Thus, we can identify the
states of $A$ with states of $B$ and the same argument leads us to the
commutative diagram of $\{\iota\}\ddownarrow \K$
\[
\xymatrix
{
\{\iota\} \sqcup B^0 \fD{} \fR{} && \{\iota\} \sqcup B \fR{} \fD{} && Y \\ 
&& && \\
\bigsqcup\limits_{\beta\in B^0} P_\beta \fR{} && \cocartesien \widehat{B} \ar@{-->}[rruu]^-{\exists !} &&
}
\]
where $B^0$ is the set of states of $B$ and where the sum
$\bigsqcup_{\beta\in B^0}$ is taken in the category
$\{\iota\}\ddownarrow \K$ (it is understood that we choose for
$\beta\in A^0$ the same $P_\alpha$ as above).  We obtain the
commutative diagram of $\{\iota\}\ddownarrow \K$:
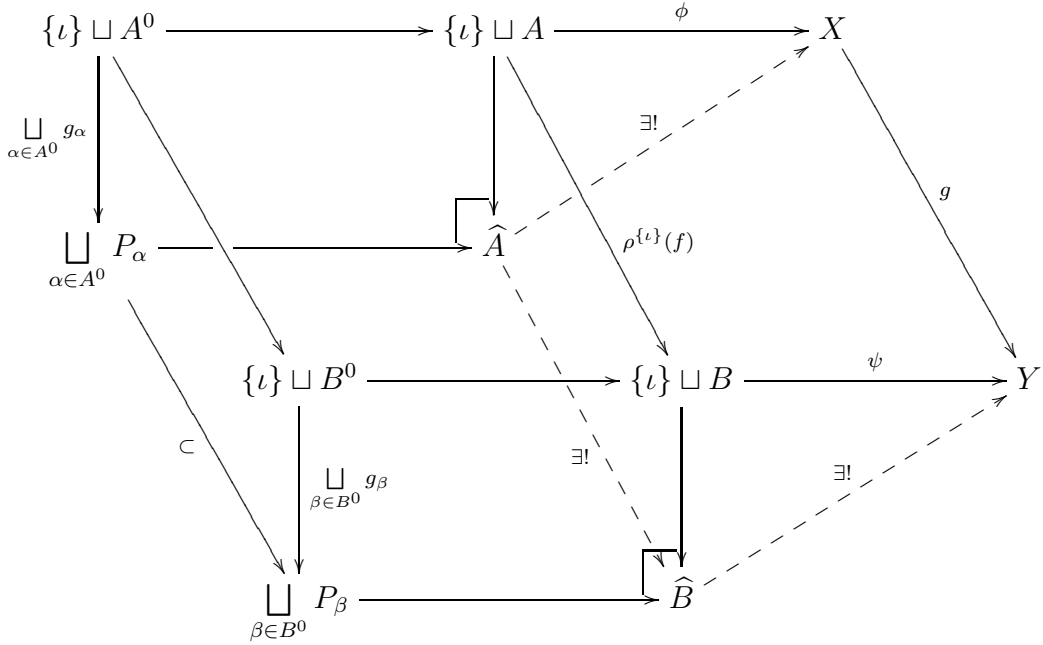
\[
\xymatrix
{
\{\iota\} \sqcup A^0 \ar@{->}[rddd]\ar@{->}[dd]_-{\bigsqcup\limits_{\alpha\in A^0} g_\alpha} \ar@{->}[rr]^-{} && \{\iota\} \sqcup A \ar@{->}[rddd]^/20pt/{\rho^{\{\iota\}}(f)} \ar@{->}[rr]^-{\phi} \ar@{->}[dd]_-{} && X \ar@{->}[rdrdd]^-{g}  \\ 
&& &&  \\
\bigsqcup\limits_{\alpha\in A^0} P_\alpha \ar@{->}[rddd]_-{\subset} \ar@{->}[rr]^-{}|(0.315)\hole && \cocartesien
\widehat{A}\ar@{-->}[rddd]_/10pt/{\exists !} \ar@{-->}[rruu]^-{\exists !}|(0.26)\hole && \\
&\{\iota\} \sqcup B^0 \ar@{->}[dd]^-{\bigsqcup\limits_{\beta\in B^0} g_\beta} \ar@{->}[rr]^-{} && \{\iota\} \sqcup B \ar@{->}[rrr]^-{\psi} \ar@{->}[dd]_-{} &&& Y \\ 
&&& && \\
&\bigsqcup\limits_{\beta\in B^0} P_\beta \ar@{->}[rr]^-{} && \cocartesien
\widehat{B} \ar@{-->}[rrruu]^-{\exists !} &&
}
\]
The map $\widehat{A} \to \widehat{B}$ making the diagram commutative
exists by the universal property of the pushout and it is one-to-one
on states and on actions, i.e. it is a cofibration. We obtain the
factorization
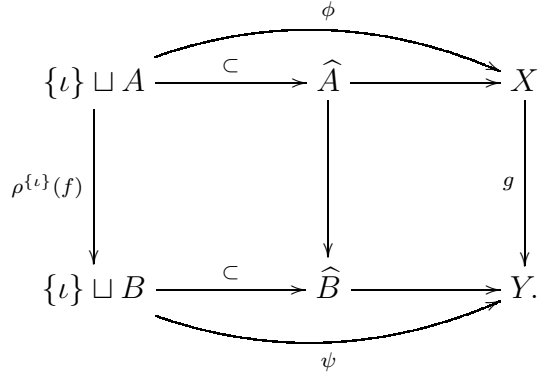
\[
\xymatrix
{\{\iota\}\sqcup A \ar@/^20pt/@{->}[rrrr]^-{\phi} \fR{\subset} \fD{\rho^{\{\iota\}}(f)} && \widehat{A} \fD{}\fR{} && X \fD{g} \\
&& &&\\
\{\iota\}\sqcup B \ar@/_20pt/@{->}[rrrr]_-{\psi} \fR{\subset} && \widehat{B} \fR{} && Y.
}
\]
By construction, the transition systems $\widehat{A}$ and
$\widehat{B}$ are star-shaped. The map of star-shaped transition
systems $\widehat{A}\to \widehat{B}$ is obtained by choosing for each
state of $B$ a map of the set $\{g:\{\iota,\alpha\} \to P\} \cup
\{\{\iota,\alpha\}\to \{\iota\}\}$. We have therefore constructed a
solution set of cofibrations for the set of generating cofibrations of
$\{\iota\}\ddownarrow \K$ with respect to $\K_{\bullet}$, i.e.  there
exists a set $J$ of cofibrations of $\{\iota\}\ddownarrow \K$ between
star-shaped objects such that every map $i\to w$ from a generating
cofibration $i$ of $\{\iota\}\ddownarrow \K$ to a map $w$ of
$\K_\bullet$ factors as a composite $i\to j\to w$ with $j\in J$.  The
proof is complete thanks to \cite[Lemma~A.3]{cubicalhdts}.  \epf

\bth \label{star-shaped} Let $\K$ be $\wts$ or $\cts$ or $\rts$. There
exists a left determined Olschok model structure on the category
$\K_{\bullet}$ of star-shaped weak (cubical, regular resp.)
transition systems with respect to the class of maps such that the
underlying map is a cofibration of $\K$.  \eth

Note that unlike in the proof of \cite[Theorem~6.8]{leftdet}, we
cannot use \cite[Theorem~4.3]{leftdet}. Indeed, $\K_\bullet$ is still
a small cone-injectivity class by Proposition~\ref{cone}. However, the
cone contains the map $\{\iota,\alpha\}\to \iota$ which is not a
cofibration in this paper.

\bpf Thanks to Theorem~\ref{comma-cat},
Proposition~\ref{corefl_bullet}, Corollary~\ref{locpre_bullet} and
Theorem~\ref{cofibrantly-gen}, we can apply
\cite[Theorem~4.1]{leftdet} if we can prove that the cylinder functor
$C_{\{\iota\}}:\{\iota\} \ddownarrow \K \to \{\iota\} \ddownarrow \K$
of $\{\iota\} \ddownarrow \K$ takes a star-shaped (weak, cubical,
regular resp.)  transition system to a star-shaped one.  Let
$\{\iota\}\to X$ be an object of $\K_{\bullet}$. We have the pushout
diagram of $\K$:
\[
\xymatrix {
  \{\iota\}\sqcup \{\iota\} \fR{} \fD{} && \{\iota\} \fD{} \\
  && \\
  C(X) \fR{} && \cocartesien \omega^{\{\iota\}}(C_{\{\iota\}}(\{\iota\}
  \to X))  }
\]
where $C$ denotes the cylinder of $\K$.  Therefore if a state $\alpha$
is reachable from $\iota$, then the state $(\alpha,\epsilon)$ with
$\epsilon=0,1$ is reachable from
$(\iota,\epsilon)=(\iota,0)=(\iota,1)$ in $C_{\{\iota\}}(\iota \to
X)$.  \epf

Let us now reconsider the argument of Section~\ref{boum}. We obtain
what follows (the functor $P_{\{\iota\}}:\{\iota\} \ddownarrow \K \to
\{\iota\} \ddownarrow \K$ denoting the right adjoint to the functor
$C_{\{\iota\}}$). Let ${\{\iota\}}\to X$ be a star-shaped transition
system of $\K_\bullet$ which contains at least one transition. Let
$(\alpha,u_1,\dots,u_n,\beta)$ be a transition of a fibrant
replacement $({\{\iota\}}\to X)^{fib}$ of ${\{\iota\}}\to X$ in
$\K_\bullet$. Let $(\gamma,\delta)$ be a pair of states of
$({\{\iota\}}\to X)^{fib}$.  If $(\alpha,\gamma)$ and $(\beta,\delta)$
are two reachable states of $P_{\{\iota\}}({\{\iota\}}\to X)$, then
the triple $(\gamma,u_1,\dots,u_n,\delta)$ is a transition of
$({\{\iota\}}\to X)^{fib}$. The crucial difference with
Section~\ref{boum} is that $(\alpha,\gamma)$ and $(\beta,\delta)$ must
now be reachable states of $P_{\{\iota\}}(({\{\iota\}}\to X)^{fib})$,
and \emph{not} any pair of states of $X^{fib}$. We have to understand
now the intuitive meaning of a reachable state of
$P_{\{\iota\}}(({\{\iota\}}\to X)^{fib})$.

Let $(\kappa,\lambda)$ be a reachable state of
$P_{\{\iota\}}(({\{\iota\}}\to X)^{fib})$. That means that there
exists a finite sequence of transitions $t_i$ of $P(X)$ (the path
space of $X$ in $\K$) from $(\alpha_i,\alpha'_i)$ to
$(\alpha_{i+1},\alpha'_{i+1})$ for $0\leq i \leq n$ with $n\geq 0$,
with $(\alpha_0,\alpha'_0) = (\iota,\iota)$ and
$(\alpha_{n+1},\alpha'_{n+1}) = (\kappa,\lambda)$. By definition of a
transition in $P(X^{fib})$, that means not only that the states
$\kappa$ and $\lambda$ are reachable, but also that the transitions
relating ${\iota}$ to $\kappa$ have the same labels as the
transitions relating ${\iota}$ to $\lambda$. Indeed, by
Proposition~\ref{cocyl-def} for $\wts$, by Proposition~\ref{inclusion}
for $\cts$ and by Proposition~\ref{path-regular} for $\rts$, the set
of actions of the path space of $X^{fib}$ is a subset of $L\p_\Sigma
L$ where $L$ is the set of actions of $X^{fib}$. Roughly speaking,
\emph{the states $\kappa$ and $\lambda$ have the same past}.

The interaction between the fibrant replacement in $\K_\bullet$ and
the causal structure can now be reformulated in plain English as
follows:

\emph{Let $(\alpha,u_1,\dots,u_n,\beta)$ be a transition of a fibrant
  replacement $({\{\iota\}}\to X)^{fib}$ of ${\{\iota\}}\to X$ in $\K_\bullet$. Let
  $(\gamma,\delta)$ be a pair of states of $({\{\iota\}}\to X)^{fib}$ such that
  $\alpha$ and $\gamma$ ($\beta$ and $\delta$ resp.) have the same
  past. Then the triple $(\gamma,u_1,\dots,u_n,\delta)$ is a
  transition of $({\{\iota\}}\to X)^{fib}$.}

\appendix

\section{Relocating maps in a transfinite composition} \label{relocation}

For this section, $\K$ is a locally presentable category and
$R$ is a map such that all maps of
$\cell_\K(\{R\})$ are epic. Proposition~\ref{dec-unique}
and Theorem~\ref{move_bad} are used as follows:
\begin{itemize}
\item In the proof of Proposition~\ref{small-mono} with the category
  $\wts$ and with the map $R:\{0,1\}\to \{0\}$.
\item In the proof of Theorem~\ref{explicit} with the category $\cts$
  and with the map $R:\{0,1\}\to \{0\}$.
\end{itemize}

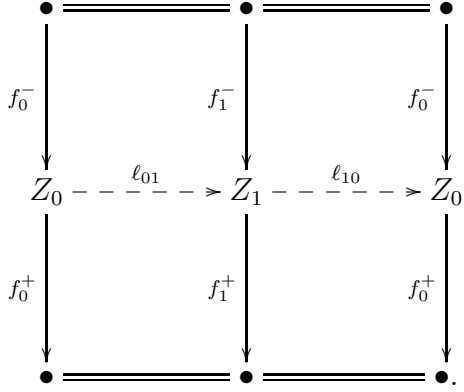
\begin{figure}
\[
\xymatrix
{
\bullet \ar@{=}[rr] \fD{f_0^-} && \bullet \fD{f_1^-} \ar@{=}[rr]&& \bullet \fD{f_0^-}\\
  &&  \\
Z_0  \ar@{-->}[rr]^-{\ell_{01}}\fD{f_0^+} && Z_1 \ar@{-->}[rr]^-{\ell_{10}}\fD{f_1^+} && Z_0\fD{f_0^+}\\
  &&  && \\
\bullet \ar@{=}[rr]  && \bullet \ar@{=}[rr]&& \bullet.
}
\]
\caption{Unique factorization $f=f^+f^-$.}
\label{unique}
\end{figure}

\bp \label{dec-unique} Every map $f$ of $\K$ factors functorially as a
composite $f=f^+f^-$ with $f^-\in \cell_\K(\{R\})$ and $f^+
\in \inj_\K(\{R\})$. This factorization is unique up to
isomorphism.  \ep

\bpf The existence of the functorial factorization is a consequence of
\cite[Proposition~1.3]{MR1780498}.  Consider the commutative diagram
of solid arrows of $\K$ of Figure~\ref{unique} with
$f^+_0f^-_0=f^+_1f^-_1$.  The lift $\ell_{01}$ exists since $f_0^-\in
\cell_\K(\{R\})$ and $f_1^+ \in
\inj_{\K}(\{R\})$. If $\ell_{01}'$ is another lift, then
there is the equality $\ell_{01} f_0^-=f_1^-= \ell_{01}' f_0^-$. By
hypothesis, the map $f_0^-$ is epic.  Therefore
$\ell_{01}=\ell_{01}'$, which means that the lift $\ell_{01}$ is
unique. By switching the two columns, we obtain another lift
$\ell_{10}$. By uniqueness of the lift, the composite
$\ell_{10}\ell_{01}$ is equal to the lift $\ell_{00}=\id_{Z_0}$ and
the composite $\ell_{01}\ell_{10}$ is equal to the lift
$\ell_{11}=\id_{Z_1}$.  \epf

\begin{figure}
\[
\xymatrix
{
P \fR{\phi} \ar@{->}[dd]^-{f^-} \ar@/_40pt/@{->}[dddd]^-{f} && X_\alpha \ar@/^40pt/@{->}[dddddd]^-/-20pt/{p_\alpha=p_{\alpha+1}q_\alpha} \ar@{->}[rrrr]^-{p_\alpha^-} \ar@{->}[dddd]^-/20pt/{q_\alpha} && && \overline{X}_\alpha \ar@{->}[dddd]^-{\overline{q}_\alpha} \ar@/^40pt/@{->}[dddddd]^-{p_\alpha^+} \\
&& && && \\
\bullet \ar@{->}[dd]^-{f^+} \ar@{-->}[rrrrrruu]^/-70pt/{\ell} && && && \\
&& && && \\
Q \fR{\psi} && \cocartesien X_{\alpha+1} \ar@{->}[rrrr]^-{\overline{\psi}} \fD{p_{\alpha+1}} && && \cocartesien \overline{X}_{\alpha+1} \fD{}\\
&& && && \\
&& \liminj X_\alpha \ar@{=}[rrrr] && && \liminj X_\alpha.
}
\]
\caption{Modification of the $\alpha$-th stage of the transfinite tower.}
\label{alpha}
\end{figure}

\bth \label{move_bad} With the notations of
Proposition~\ref{dec-unique}. Let $A$ be a set of maps of $\K$ such
that $A\subset \inj_\K(\{R\})$.  Then every map $f\in
\cell_\K(A \cup \{R\})$ factors uniquely, up to
isomomorphism, as a composite $f=f^+f^-$ with $f^-\in
\cell_\K(\{R\})$ and $f^+\in \cell_\K(A)$.  \eth

Theorem~\ref{move_bad} means that the cells $R$ of a
cellular complex $\cell_\K(A \cup \{R\})$ can be relocated
at the beginning of the cellular decomposition.

\bpf Let $(q_\alpha:X_\alpha\to X_{\alpha+1})_{\alpha \geq 0}$ be a
transfinite tower of pushouts of maps of $A \cup \{R\}$.
Consider the commutative diagram of solid arrows of
Figure~\ref{alpha}. It represents the $\alpha$-th stage of the tower
$(q_\alpha:X_\alpha \to X_{\alpha+1})_{\alpha \geq 0}$ which is
supposed to be a pushout of a map $f$ of $A \cup \{R\}$.
Since the factorizations $f=f^+f^-$ and
$p_\alpha=p_\alpha^+p_\alpha^-$ are functorial, there exists a map
$\ell:\bullet \to \overline{X}_\alpha$ such that $\ell f^-=p^-_\alpha
\phi$.  We obtains $\overline{q}_\alpha\ell f^-=\overline{q}_\alpha
p^-_\alpha \phi = \overline{\psi}\psi f= \overline{\psi}\psi f^+ f^-$.
By hypothesis, the map $f^-$ is epic. We obtain
$\overline{q}_\alpha\ell = \overline{\psi}\psi f^+$.  Finally, an
immediate application of the universal property of a pushout square
shows that the commutative square
\[
\xymatrix {
  \bullet \fR{\ell}\ar@{->}[dd]^-{f^+} && \overline{X}_\alpha \ar@{->}[dd]^-{\overline{q}_\alpha}\\
&& \\
Q\fR{\overline{\psi}\psi} && \cocartesien \overline{X}_{\alpha+1}} 
\]
is a pushout square. This process can be iterated by composing the
$(\alpha+1)$-th attaching map with the map
$\overline{\psi}:X_{\alpha+1} \to \overline{X}_{\alpha+1}$.  We have
obtained by induction on $\alpha \geq 0$ a new tower
$(\overline{q}_\alpha:\overline{X}_\alpha\to
\overline{X}_{\alpha+1})_{\alpha \geq 0})$ with the same colimit and a
map of towers $q_* \to \overline{q}_*$. Consequently, the map $p_0:X_0
\to \liminj X_\alpha$ factors as a composite \[p_0 : X_0
\stackrel{p_0^-} \longrightarrow \overline{X}_0 \stackrel{p_0^+}
\longrightarrow \liminj X_\alpha\] such that the right-hand map is a
transfinite composition of pushouts of maps of the set $\{f^+ \mid
f\in A \cup \{R\} \}$. There is the equality
$R^-=R$ and therefore $R^+=\id$.
By hypothesis, there is the inclusion $A\subset
\inj_\K(\{R\})$, which implies $f=f^+$ for all $f\in
A$. Thus, there is the equality $\{f^+ \mid f\in A \cup
\{R\}\} = A \cup \{\id\}$ and the proof is complete.  \epf


\begin{thebibliography}{Gau15b}

\bibitem[AHS06]{topologicalcat}
J.~Ad{\'a}mek, H.~Herrlich, and G.~E. Strecker.
\newblock Abstract and concrete categories: the joy of cats.
\newblock {\em Repr. Theory Appl. Categ.}, (17):1--507 (electronic), 2006.
\newblock Reprint of the 1990 original [Wiley, New York; MR1051419].

\bibitem[AR94]{MR95j:18001}
J.~Ad{\'a}mek and J.~Rosick{\'y}.
\newblock {\em Locally presentable and accessible categories}.
\newblock Cambridge University Press, Cambridge, 1994.

\bibitem[Bek00]{MR1780498}
T.~Beke.
\newblock Sheafifiable homotopy model categories.
\newblock {\em Math. Proc. Cambridge Philos. Soc.}, 129(3):447--475, 2000.

\bibitem[CS96]{MR1461821}
G.~L. Cattani and V.~Sassone.
\newblock Higher-dimensional transition systems.
\newblock In {\em 11th Annual IEEE Symposium on Logic in Computer Science (New
  Brunswick, NJ, 1996)}, pages 55--62. IEEE Comput. Soc. Press, Los Alamitos,
  CA, 1996.

\bibitem[Gau08]{ccsprecub}
P.~Gaucher.
\newblock Towards a homotopy theory of process algebra.
\newblock {\em Homology Homotopy Appl.}, 10(1):353--388 (electronic), 2008.

\bibitem[Gau10]{hdts}
P.~Gaucher.
\newblock Directed algebraic topology and higher dimensional transition
  systems.
\newblock {\em New York J. Math.}, 16:409--461 (electronic), 2010.

\bibitem[Gau11]{cubicalhdts}
P.~Gaucher.
\newblock Towards a homotopy theory of higher dimensional transition systems.
\newblock {\em Theory Appl. Categ.}, 25:No.\ 25, 295--341 (electronic), 2011.

\bibitem[Gau14]{homotopyprecubical}
P.~Gaucher.
\newblock Homotopy theory of labelled symmetric precubical sets.
\newblock {\em New York J. Math.}, 20:93--131 (electronic), 2014.

\bibitem[Gau15a]{csts}
P.~Gaucher.
\newblock The geometry of cubical and regular transition systems.
\newblock {\em Cah. Topol. G\'eom. Diff\'er. Cat\'eg.}, LVI-4, 2015.

\bibitem[Gau15b]{leftdet}
P.~Gaucher.
\newblock Left determined model categories.
\newblock arXiv:1507.02128, to appear in New York J. of Math.

\bibitem[Hir03]{ref_model2}
P.~S. Hirschhorn.
\newblock {\em Model categories and their localizations}, volume~99 of {\em
  Mathematical Surveys and Monographs}.
\newblock American Mathematical Society, Providence, RI, 2003.

\bibitem[Hir15]{undercat}
P.~Hirschhorn.
\newblock Overcategories and undercategories of model categories.
\newblock arXiv:1507.01624, 2015.

\bibitem[Hov99]{MR99h:55031}
M.~Hovey.
\newblock {\em Model categories}.
\newblock American Mathematical Society, Providence, RI, 1999.

\bibitem[KR05]{ideeloc}
A.~Kurz and J.~Rosick{\'y}.
\newblock Weak factorizations, fractions and homotopies.
\newblock {\em Applied Categorical Structures}, 13(2):pp.141--160, 2005.

\bibitem[Ols09a]{MO}
M.~Olschok.
\newblock Left determined model structures for locally presentable categories.
\newblock {\em Applied Categorical Structures}, 2009.
\newblock 10.1007/s10485-009-9207-2.

\bibitem[Ols09b]{MOPHD}
M.~Olschok.
\newblock {\em On constructions of left determined model structures}.
\newblock PhD thesis, Masaryk University, Faculty of Science, 2009.

\bibitem[Ros09]{MR2506258}
J.~Rosick{\'y}.
\newblock On combinatorial model categories.
\newblock {\em Appl. Categ. Structures}, 17(3):303--316, 2009.

\end{thebibliography}
\end{document}